\documentclass[10pt]{article}

\usepackage{amsmath}
\usepackage{amsfonts}
\usepackage{amssymb}
\usepackage[english]{babel}
\usepackage{amsthm}
\usepackage{extarrows}
\usepackage{amscd}
\usepackage{tikz}
\usepackage{mathdots}
\usetikzlibrary{arrows,shapes}

\def\pf{\par\noindent {\bf Proof}~\par\noindent}

\newcommand{\mR}{\mathbb{R}}

\newcommand{\mN}{\mathbb{N}}

\newcommand{\mE}{\mathbb{E}}
\newcommand{\mS}{\mathbb{S}}

\newcommand{\mT}{\mathbb{T}}
\newcommand{\mU}{\mathbb{U}}

\newcommand{\mcD}{\mathcal{D}}

\newcommand{\mcV}{\mathcal{V}}
\newcommand{\mcW}{\mathcal{W}}

\newcommand{\ux}{\underline{x}}
\newcommand{\uy}{\underline{y}}

\newcommand{\uom}{\underline{\omega}}

\newcommand{\uD}{\underline{D}}
\newcommand{\uS}{\underline{S}}

\newcommand{\p}{\partial}
\newcommand{\dirac}{\underline{\p}}

\newcommand{\eop}{\hfill$\square$}

\newcommand{\onehalf}{\frac{1}{2}}

\newcommand{\invr}{\frac{1}{r}}
\newcommand{\invrsq}{\frac{1}{r^2}}
\newcommand{\invux}{\frac{1}{\ux}}

\newtheorem{lemma}{Lemma}
\newtheorem{proposition}{Proposition}
\newtheorem{definition}{Definition}
\newtheorem{remark}{Remark}

\newtheorem{example}{Example}

\begin{document}

\title{Distributions: a spherical co-ordinates approach}
\author{Fred Brackx}

\date{Clifford Research Group, Department of Mathematical Analysis,\\
Faculty of Engineering and Architecture, Ghent University
}

\maketitle


\begin{abstract}
When expressing a distribution in Euclidean space in spherical co-ordinates, derivation with respect to the radial and angular co-ordinates is far from trivial. Exploring the possibilities of defining a radial derivative of the delta-distribution $\delta(\ux)$ (the angular derivatives of $\delta(\ux)$ being zero since the delta-distribution is itself radial) led, see \cite{bsv}, to the introduction of a new kind of distributions, the so--called {\em signumdistributions}, as continuous linear functionals on a space of test functions showing a singularity at the origin. In this paper we search for a definition of  the radial and angular derivatives of a general standard distribution and  again, as expected, we are inevitably led to consider signumdistributions. Although these signumdistributions provide an adequate framework for the actions on distributions aimed at, it turns out that the derivation with respect to the radial distance of a general (signum)distribution is still not yet unique, but gives rise to an equivalent class of (signum)distributions.\\

\end{abstract}


\section{Introduction}
\label{intro}


Let us consider a scalar--valued distribution $T(\ux) \in \mcD'(\mR^m)$ expressed in terms of spherical co-ordinates: $\ux = r \uom,\, r = |\ux|,\, \uom = \sum_{j=1}^{m} \, e_j \, \omega_j \, \in \mS^{m-1}$, $(e_j)_{j=1}^m$ being an orthonormal basis of $\mR^m$ and $\mS^{m-1}$ being the unit sphere in $\mR^m$. The aim of this paper is to search for an adequate definition of the radial and angular derivatives $\p_r\, T$ and $\p_{\omega_j} \, T, \, j=1,\ldots,m$. This problem was treated in \cite{bsv} for the special and interesting case of the delta-distribution $\delta(\ux)$, the following spherical co-ordinates expression of which is often encountered in physics texts:
\begin{equation}
\label{physicsdelta}
\delta(\ux) = \frac{1}{a_m}Ê\, \frac{\delta(r)}{r^{m-1}}
\end{equation}
with $a_m = \frac{2\pi^{\frac{m}{2}}}{\Gamma(\frac{m}{2})}$ the area of the unit sphere $\mS^{m-1}$ in $\mR^m$.
Apparently this expression (\ref{physicsdelta}) can be explained in the following way. Write the action of the delta-distribution as an integral:
\begin{eqnarray*}
\varphi(0) = \langle \ \delta(\ux) , \varphi(\ux) \ \rangle &=& \int_{\mR^m} \, \delta(\ux) \, \varphi(\ux) \, dV(\ux)\\
&=& \int_{0}^{\infty} \, r^{m-1} \delta(\ux) \, dr \, \int_{\mS^{m-1}} \, \varphi(r \, \uom) \, dS_{\uom}\\
&=& a_m\, \int_{0}^\infty \, r^{m-1} \, \delta(\ux) \, \Sigma^{0}[\varphi](r) \, dr
\end{eqnarray*}
introducing the so--called {\em spherical mean} of the test function $\varphi$ given by
$$
\Sigma^{0}[\varphi](r) = \frac{1}{a_m}Ê\, \int_{\mS^{m-1}} \, \varphi(r \, \uom) \, dS_{\uom}
$$
As it is easily seen that $\Sigma^0[\varphi](0) = \varphi(0)$ it follows that
$$
a_m\, \int_{0}^\infty \, r^{m-1} \, \delta(\ux) \, \Sigma^{0}[\varphi](r) \, dr = \int_{0}^\infty \, \delta(r) \, \Sigma^{0}[\varphi](r) \, dr  = \langle \  \delta(r) ,  \Sigma^{0}[\varphi](r) \ \rangle 
$$
which explains (\ref{physicsdelta}).
However we prefer to interpret this expression (\ref{physicsdelta}) mathematically as
\begin{equation}
\label{sigmanul}
\varphi(0) = \langle \ \delta(\ux) , \varphi(\ux) \ \rangle = \langle \  \delta(r) ,  \Sigma^{0}[\varphi](r) \ \rangle = \Sigma^{0}[\varphi](0)
\end{equation}
Straightforward successive derivation with respect to $r$ of (\ref{physicsdelta}) leads to
\begin{eqnarray}
\label{evendelta}
\p_r^{2\ell}\,\delta(\ux) &=& \frac{1}{(2\ell)!}(m)(m+1)\cdots(m+2\ell-1)\, \frac{1}{a_m}\, \frac{\delta^{(2\ell)}(r)}{r^{m-1}}\\
\label{odddelta}
\p_r^{2\ell+1}\,\delta(\ux) &=& \frac{1}{(2\ell+1)!}(m)(m+1)\cdots(m+2\ell)\, \frac{1}{a_m}\, \frac{\delta^{(2\ell+1)}(r)}{r^{m-1}}
\end{eqnarray}
Expression (\ref{evendelta}) then is interpreted as
$$
\langle \ \p_r^{2\ell}\,\delta(\ux) ,  \varphi(\ux) \ \rangle = \frac{1}{(2\ell)!}(m)(m+1)\cdots(m+2\ell-1)\, \langle \  \delta^{(2\ell)}(r) , \Sigma^{0}[\varphi](r) \ \rangle
$$
which is meaningful and which can serve as the definition of the even order derivatives with respect to $r$ of the delta-distribution.
However expression (\ref{odddelta}) makes no sense since the spherical mean $\Sigma^{0}[\varphi](r)$ is an even function of $r$, whence its odd order derivatives vanish at the origin:
$$
 \langle \  -\,\p_r^{2\ell+1}\, \delta(r) ,  \Sigma^{0}[\varphi](r) \ \rangle =  \{\p_r^{2\ell+1}\,\Sigma^{0}[\varphi](r)\}|_{r=0} = 0
$$
So the question arises how to explain the fact that, proceeding stepwise by derivation with respect to $r$, the even order derivatives of $\delta(\ux)$ apparently make sense, while its odd order derivatives are zero distributions, in this way violating the basic requirement of any derivation procedure that $\p_r\, \p_r$ should equal $\p_r^2$. Let us to that end have a quick look at the functional analytic background of this phenomenon; for a more systematic treatment we refer to \cite{bsv}.\\

When expressing a scalar--valued test function $\varphi(\ux) \in \mcD(\mR^m)$ in spherical co-ordinates, one obtains a function $\widetilde{\varphi}(r, \uom) = \varphi(r\uom) \in 
\mcD(\mR \times \mS^{m-1})$, but it is evident that not all functions $\widetilde{\varphi}(r, \uom)  \in \mcD(\mR \times \mS^{m-1})$ stem from a test function in $\mcD(\mR^m)$.
However a one--to--one correspondence may be established between the usual space of test functions $\mcD(\mR^m)$ and a specific subspace of $\mcD(\mR \times \mS^{m-1})$.

\begin{lemma} (see \cite{helgason})
\label{iso}
There is a one--to--one correspondence $\varphi(\ux) \leftrightarrow \widetilde{\varphi}(r, \uom) = \varphi(r\uom)$ between the spaces $\mcD(\mR^m)$ and $\mcV = \{\phi(r,\uom) \in \mcD(\mR \times \mS^{m-1})  : \phi$ is even, i.e. $\phi(-r,-\uom) = \phi(r,\uom)$, and  $\{\p_r^n \, \phi(r,\uom) \}|_{r=0}$ is a homogeneous polynomial of degree
$n$ in $(\omega_1,\ldots,\omega_m), \forall n \in \mN\}$.
\end{lemma}

\noindent
Clearly $\mcV$ is a closed (but not dense) subspace of $\mcD(\mR \times \mS^{m-1})$ and even of $\mcD_E(\mR \times \mS^{m-1})$, where the subscript $E$ refers to the even character of the test functions in that space; this space $\mcV$ is endowed with the induced topology of $\mcD(\mR \times \mS^{m-1})$. The one--to--one correspondence between the spaces of test functions $\mcD(\mR^m)$ and $\mcV$ translates into a one--to--one correspondence between the standard distributions $T \in \mcD'(\mR^m)$ and the bounded linear functionals in $\mcV'$, this correspondence being given  by
$$
\langle \ T(\ux) , \varphi(\ux) \ \rangle = \langle \ \widetilde{T}(r,\uom) , \widetilde{\varphi}(r,\uom) \ \rangle
$$
By Hahn--Banach's theorem the bounded linear functional $\widetilde{T}(r,\uom)  \in \mcV'$ may be extended to the distribution $\mT(r,\uom) \in \mcD'(\mR \times \mS^{m-1})$; such an extension is called a {\em spherical representation} of the distribution $T$ (see e.g.\ \cite{vindas}). However as the subspace $\mcV$ is not dense in  $\mcD(\mR \times \mS^{m-1})$, the spherical representation of a distribution is {\em not unique}, but if $\mT_1$ and $\mT_2$ are two different spherical representations of the same distribution $T$, their restrictions to $\mcV$ coincide: 
$$
\langle \ \mT_1(r,\uom) , \widetilde{\varphi}(r,\uom) \ \rangle = \langle \ \mT_2(r,\uom) , \widetilde{\varphi}(r,\uom) \ \rangle = 
\langle \ \widetilde{T}(r,\uom) , \varphi(r\uom) \ \rangle = \langle \ T(\ux) , \varphi(\ux) \ \rangle
$$

\noindent
For test functions in $\mcD(\mR \times \mS^{m-1})$ the spherical variables $r$ and $\uom$ are ordinary variables, and thus smooth functions. It follows that for distributions in $\mcD'(\mR \times \mS^{m-1})$ multiplication by $r$ and $\omega_j, j=1,\ldots,m$ and differentiation with respect to $r$ and $\omega_j, j=1,\ldots,m$ are well--defined standard  operations, whence
\begin{equation}
\label{blackboardaction}
\langle \ \p_r\, \mT(r,\uom) , \Xi(r,\uom) \ \rangle = -\, \langle \ \mT(r,\uom) , \p_r\,  \Xi(r,\uom) \ \rangle 
\end{equation}
for all test functions $\Xi(r,\uom) \in \mcD(\mR \times \mS^{m-1})$, and similar expressions for $\p_{\omega_j}\, \mT$, $r\, \mT$ and $\uom\, \mT$.
However if $\mT_1$ and $\mT_2$ are two different spherical representations of the same distribution $T \in \mcD'(\mR^m)$, then, upon restriction to test functions $\widetilde{\varphi}(r,\uom) \in \mcV$, we are stuck with
$$
-\, \langle \ \mT_1(r,\uom) , \p_r\, \widetilde{\varphi}(r,\uom)  \ \rangle  \neq -\, \langle \ \mT_2(r,\uom) , \p_r\, \widetilde{\varphi}(r,\uom)  \ \rangle 
$$
because $\p_r\, \widetilde{\varphi}(r,\uom)$  does no longer belong to $\mcV$ (and neither do $\p_{\omega_j}\, \widetilde{\varphi}(r,\uom)$, $r\, \widetilde{\varphi}(r,\uom)$ and $\uom\, \widetilde{\varphi}(r,\uom)$) since it is an odd function in the variables $(r,\uom)$. And it is also clear that the action (\ref{blackboardaction}) might be unambiguously restricted to testfunctions in $\mcV$ if the test function $\Xi$ were in a subspace of $\mcD(\mR \times S^{m-1})$ consisting of odd functions.
The conclusion is that the concept of spherical representation of a distribution does not allow for an unambiguous definition of the actions proposed. What is more, it becomes apparent that there is a need for a subspace of odd test functions.
And at the same time it becomes clear why even order derivatives with respect to $r$ of the delta-distribution and of a standard distribution in general, are well-defined instead. Indeed, we have e.g. 
$$
\langle \ \p_r^{2\ell} \, \mT(r,\uom) , \Xi(r,\uom) \ \rangle =  \langle \ \mT(r,\uom) , \p_r^{2\ell} \,  \Xi(r,\uom) \ \rangle 
$$
where now $\p_r^{2\ell}\,  \Xi(r,\uom)$ belongs to $\mcD_E(\mR \times \mS^{m-1})$ which enables restriction to test functions in  $\mcV$ in an unambiguous way.


\section{Signumdistributions}
\label{signumdistributions}


As already observed in the preceding section, $\uom$ is an ordinary (vector) variable in $\mR \times \mS^{m-1}$, whence it makes sense to consider the following subspace of vector--valued test functions in $\mR \times \mS^{m-1}$:
$$
\mcW  = \uom\, \mcV \subset \mcD_O(\mR \times \mS^{m-1}; \mR^m) \subset \mcD(\mR \times \mS^{m-1};\mR^m)
$$
where now the subscript $O$ refers to the odd character of the test functions under consideration, i.e. $\psi(-r,-\uom) = -\, \psi(r,\uom), \forall \psi \in  \mcD_O(\mR \times \mS^{m-1}; \mR^m)$. This space $\mcW$ is endowed with the induced topology of $\mcD(\mR \times \mS^{m-1};\mR^m)$. By definition there is a one--to--one correspondence between the spaces $\mcV$ and $\mcW$.\\

\noindent
From now on we will interpret vectors in $\mR^m$ as Clifford $1$--vectors in the Clifford algebra $\mR_{0,m}$, where the basis vectors $(e_j, j=1,\ldots,m)$ of $\mR^m$, satisfy the relations $e_j^2 = - 1,\, e_i \wedge e_j = e_ie_j = - e_je_i = - e_j \wedge e_i,\, e_i \cdot e_j = 0, i \neq j=1,\ldots,m$. This allows for the use of the very efficient {\em geometric} or {\em Clifford product} of Clifford vectors:
$$
\ux \,  \uy = \ux \cdot \uy  + \ux \wedge \uy
$$
for which, in particular,
$$
\ux \,  \ux = \ux \cdot \ux  = -\, |\ux|^2
$$
$\ux$ being the Clifford 1--vector $\ux = \sum_{j=1}^m \, e_j\, x_j$, whence also
$$
\uom \,  \uom = \uom \cdot \uom  = -\, |\uom|^2 = -1
$$
For more on Clifford algebras we refer to e.g.\ \cite{port}.\\

\noindent
For each $\mU(r,\uom) \in \mcD'(\mR \times \mS^{m-1};\mR^m)$ we define $\widetilde{U}(r,\uom) \in \mcW'$ by the restriction
$$
\langle \ \widetilde{U}(r,\uom) , \uom\, \widetilde{\varphi}(r,\uom) \ \rangle =  \langle \ \mU(r,\uom) , \uom\, \widetilde{\varphi}(r,\uom) \ \rangle , \quad \forall \ \uom\, \widetilde{\varphi}(r,\uom) \in \mcW
$$
In $\mR^m$ we consider the space $\Omega(\mR^m; \mR^m) = \{ \uom\, \varphi(\ux) : \varphi(\ux) \in \mcD(\mR^m)  \}$. Clearly the functions in $\Omega(\mR^m; \mR^m)$ are no longer differentiable in the whole of $\mR^m$, since they are not defined at the origin due to the function $\uom = \frac{\ux}{|\ux|}$. By definition there is a one--to--one correspondence between the spaces $\mcD(\mR^m)$ and $\Omega(\mR^m; \mR^m)$.\\
For each $\widetilde{U}(r,\uom) \in \mcW'$ we define $^{s}U(\ux)$ by
$$
\langle \ ^{s}U(\ux) , \uom\, \varphi(\ux) \ \rangle =  \langle \ \widetilde{U}(r,\uom) , \uom\, \widetilde{\varphi}(r,\uom)) \ \rangle , \quad \forall \ \uom\,\varphi(\ux) \in \Omega(\mR^m; \mR^m)
$$
Clearly $^{s}U(\ux)$ is a bounded linear functional on $\Omega(\mR^m; \mR^m)$, for which, in \cite{bsv}, we coined the term {\em signumdistribution}.\\

\noindent
Now start with a standard distribution $T(\ux) \in \mcD'(\mR^m)$ and let $\mT(r,\uom) \in \mcD'(\mR \times \mS^{m-1})$ be one of its spherical representations. Put
$\mS(r,\uom) = \uom\, \mT(r,\uom)$ which in its turn leads to the signumdistribution $^{s}S(\ux) \in \Omega'(\mR^m; \mR^m)$. Then we consecutively have
\begin{eqnarray*}
\langle \  ^{s}S(\ux) ,  \uom\, \varphi(\ux)  \ \rangle = \langle \  \mS(r,\uom) , \uom\, \widetilde{\varphi}(r,\uom)    \ \rangle &=& \langle \  \uom\, \mT(r,\uom) ,  \uom\,  \widetilde{\varphi}(r,\uom)  \ \rangle\\ &=& -\,  \langle \   \mT(r,\uom) ,   \widetilde{\varphi}(r,\uom)  \ \rangle = -\, \langle \  T(\ux) ,  \varphi(\ux)  \ \rangle 
\end{eqnarray*}
since $\uom^2 = -1$. We call $^{s}S(\ux)$ a signumdistribution associated to the distribution $T(\ux)$ and denote it by $T^{\vee}(\ux)$. It thus holds that for all test functions $\varphi \in \mcD(\mR^m)$
\begin{equation}
\label{unique}
\langle \  T^{\vee}(\ux) ,  \uom\, \varphi(\ux)  \ \rangle = -\, \langle \  T(\ux) ,  \varphi(\ux)  \ \rangle 
\end{equation}
At the same time we call $T(\ux)$ the distribution associated to the signumdistribution $^{s}S(\ux)$ and we denote this distribution by $^{s}S^{\wedge}(\ux)$. Formula (\ref{unique}) then also reads
\begin{equation}
\label{unique2}
\langle \  ^{s}S(\ux) ,  \uom\, \varphi(\ux)  \ \rangle = -\, \langle \  ^{s}S^{\wedge}(\ux) ,  \varphi(\ux)  \ \rangle 
\end{equation}
and it is clear that
$$
T^{\vee \wedge} = T \qquad {\rm and} \qquad ^{s}S^{\wedge \vee}  =  {^{s}S}
$$
\noindent
At first sight for a given distribution $T(\ux)$ the associated signumdistribution $T^{\vee}(\ux)$ is not uniquely defined since its construction involves the not uniquely defined spherical representation $\mT$ of $T(\ux)$.
Nevertheless it follows from (\ref{unique}) that for a given distribution $T(\ux)$ its associated signumdistribution $T^{\vee}(\ux)$ is unique, what can also be proven directly as follows.

\begin{proposition}
Given the distribution $T(\ux)$ its associated signumdistribution $T^{\vee}(\ux)$ is uniquely determined.
\end{proposition}

\pf
Assume that $\mT_1$ and $\mT_2$ are two different spherical representations of $T$, i.e. for all test functions $\Xi(r,\uom) \in \mcD(\mR \times S^{m-1}; \mR^m)$ holds
$$
\langle \ \mT_1 , \Xi(r,\uom) \ \rangle \ \neq \ \langle \ \mT_2 , \Xi(r,\uom) \ \rangle
$$
while for all test functions $\widetilde{\varphi}(r,\uom) \in \mcV$ holds
$$
\langle \ \mT_1 , \widetilde{\varphi}(r,\uom) \ \rangle \ = \ \langle \ \mT_2 , \widetilde{\varphi}(r,\uom) \ \rangle \ = \ \langle \, \widetilde{T} , \widetilde{\varphi}(r,\uom) \ \rangle 
$$
Let $T^{\vee}_1$ and $T^{\vee}_2$ be the associated signumdistributions to $T$ through the spherical representations $\mT_1$ and $\mT_2$ respectively. Then for $j=1,2$ it holds that
$$
\langle \ T^{\vee}_j , \uom\, \varphi(\ux) \ \rangle \ = \ \langle \ \mT_j , \widetilde{\varphi}(r,\uom) \ \rangle
$$
whence $T^{\vee}_1  = T^{\vee}_2$ on $\Omega(\mR^m; \mR^m)$.
\eop\\

Conversely for a given signumdistribution $^{s}U \in \Omega'(\mR^m; \mR^m)$ we define the associated distribution $^{s}U^{\wedge}$ by
$$
\langle \  ^{s}U^{\wedge}(\ux) ,   \varphi(\ux)  \ \rangle = -\, \langle \  ^{s}U(\ux) ,  \uom\, \varphi(\ux)  \ \rangle \qquad \forall \varphi(\ux) \in \mcD(\mR^m)
$$
Clearly it holds that
$$
T^{\vee \wedge} = T \hspace{8mm} \mbox{and}Ê \hspace{8mm} ^{s}U^{\wedge \vee} = ^{s}U
$$

\noindent
As an example consider the distribution $T(\ux) =  \delta(\ux)$. Our aim is to define the signumdistribution $\delta^{\vee}(\ux)$.
A spherical representation of the delta-distribution is given by
$$
\langle \ \mT(r,\omega) , \Xi(r,\uom) \ \rangle = \Sigma^0[\Xi(r,\uom)]\}|_{r=0}
$$
Indeed, when restricting to the space $\mcV$ and taking into account property (\ref{sigmanul}), we obtain
$$
\langle \ \mT(r,\omega) , \widetilde{\varphi}(r,\uom) \ \rangle =  \Sigma^0[\varphi(r\,\uom)]\}|_{r=0} = \langle \ \delta(\ux) , \varphi(\ux) \ \rangle
$$
This particular spherical representation of $T(\ux)$ induces the signumdistribution associated to $\delta(\ux) $, which we define to be  $\delta^{\vee}(\ux)$ . It thus holds that for all test functions $\varphi \in \mcD(\mR^m)$
\begin{equation}
\label{omegadelta}
\langle \   \delta^{\vee}(\ux) ,  \uom\, \varphi(\ux)  \ \rangle =  -\, \langle \  \delta(\ux) ,  \varphi(\ux)  \ \rangle 
\end{equation}
For further examples we refer to \cite{bsv}.


\section{The Dirac operator in spherical co-ordinates}
\label{diracoperator}


The Dirac operator $\dirac = \sum_{j=1}^m \, e_j\, \p_{x_j}$, which may be seen as a Stein--Weiss projection of the gradient operator (see e.g.\ \cite{stein}) and which underlies the higher dimensional theory of the monogenic functions in so--called {\em Clifford analysis} (see e.g.\ \cite{dss, ghs}),  linearizes the Laplace operator: $\dirac^2 = -\, \Delta$. Its action  on a scalar--valued standard distribution $T(\ux)$ results into the vector--valued distribution $\dirac\, T(\ux)$ given for all test functions $\varphi(\ux) \in \mcD(\mR^m)$ by
$$
\langle \ \dirac\, T(\ux)\, , \varphi(\ux) \ \rangle = \sum_{j=1}^m \, e_j\, \langle \ \p_{x_j}\, T(\ux)\, , \varphi(\ux) \ \rangle =
- \, \sum_{j=1}^m \, e_j\, \langle \  T(\ux)\, , \p_{x_j}\, \varphi(\ux) \ \rangle = - \, \langle \  T(\ux)\, , \dirac\, \varphi(\ux) \ \rangle
$$
which is a meaningful operation since only derivatives with respect to the cartesian co-ordinates are involved.\\
Two fundamental formulae in monogenic function theory are
$$
\{ \ux , \dirac\} = \ux \, \dirac + \dirac \, \ux =  -\, 2 \mE - m \hspace{8mm} \mbox{and}Ê \hspace{8mm}  [ \ux , \dirac ] = \ux \, \dirac - \dirac \, \ux = m - 2 \Gamma
$$
where $$\mE = \sum_{j=1}^m \, x_j \, \p_{x_j}$$ is the scalar Euler operator, and $$\Gamma = \sum_{j<k} \, e_j e_k\, L_{jk} =  \sum_{j<k} \, e_j e_k (x_j \p{x_k} - x_k \p_{x_j})$$ is the bivector angular momentum operator. It follows that
$$
\ux\, \dirac = - \mE - \Gamma
$$
or more precisely
$$
 \ux  \cdot \dirac  =  - \mE \hspace{8mm}  \mbox{and}Ê\hspace{8mm} \ux \wedge \dirac  = - \Gamma
$$
Passing to spherical co-ordinates $\ux = r \uom,\, r = |\ux|,\, \uom = \sum_{j=1}^{m} \, e_j \, \omega_j \, \in \mS^{m-1}$, the Dirac operator takes the form
$$
\dirac = \dirac_{rad} + \dirac_{ang}
$$
with
$$
\dirac_{rad} = \uom\, \p_r   \hspace{8mm} \mbox{and}  \hspace{8mm}  \dirac_{ang} = \frac{1}{r}\, \p_{\uom}
$$
To give an idea what the angular differential operator $\p_{\uom} = \sum_{j=1}^m\, e_j\, \p_{\omega_j}$ looks like, let us mention its explicit form in dimension $m=2: \p_{\uom} = e_{\theta} \, \p_{\theta}$ and in dimension $m=3: \p_{\uom} = e_{\theta} \, \p_{\theta} + e_{\varphi} \, \frac{1}{\sin{\theta}}\, \p_{\varphi}$, the meaning of the polar co-ordinates $\theta$ and $\varphi$ being straightforward. The operator $\p_{\uom}$ is sometimes called the {\em spherical Dirac operator}.\\

\noindent
Taking into account that $\p_{\uom}$ is orthogonal to $\uom$, the Euler operator in spherical co-ordinates then reads: 
$$
\mE = -\,  \ux  \cdot \dirac   =  -\, r\uom \cdot \dirac_{rad} = -\, r\uom \cdot \uom\, \p_r = r \, \p_r
$$
while the angular momentum operator $\Gamma$ takes the form 
$$
\Gamma = -\,  \ux  \wedge \dirac  =  -\, r\uom \wedge \dirac_{ang} = -\, r\uom \wedge  \frac{1}{r}\, \p_{\uom} = -\, \uom \wedge \p_{\uom} = - \uom \, \p_{\uom} 
$$
The question now is how to define, if possible, the action of the $\dirac_{rad}$ and $\dirac_{ang}$ operators on a standard distribution. To that end both operators should be expressed in terms of cartesian derivatives, which is achieved by putting
$$
\dirac_{rad} = \uom\, \p_r = -\, \frac{1}{\ux} \, \mE \hspace{8mm} \mbox{and} \hspace{8mm} \dirac_{ang} = \frac{1}{r}\, \p_{\uom} = -\, \frac{1}{\ux} \, \Gamma
$$
It immediately becomes clear that, in this way, the actions of $\dirac_{rad}$ and $\dirac_{ang}$ on a standard distribution $T(\ux)$ are well-defined but not uniquely defined. Indeed, due to the division by the analytic function $\ux$ showing a zero at the origin, both expressions
\begin{equation}
\label{equirad}
\dirac_{rad}\, T(\ux) = \uom\, \p_r\, T(r \uom) =  -\, \left[\frac{1}{\ux} \, \mE\, T(\ux)\right] 
\end{equation}
and
\begin{equation}
\label{equiang}
\dirac_{ang}\, T(\ux) =  \frac{1}{r}\, \p_{\uom}\, T(r\uom) = -\, \left[ \frac{1}{\ux} \, \Gamma\, T(\ux)\right]
\end{equation}
represent equivalent classes of distributions each two of which differ by a vector multiple of the delta-distribution $\delta(\ux)$. However if $\uS_1 = \dirac_{rad}\,T(\ux)$ and $\uS_2 = \dirac_{ang}\, T(\ux)$ are distributions arbitrarily chosen in the equivalent classes (\ref{equirad}) and (\ref{equiang})  respectively, i.e.
$$
\ux \, \uS_1 = -\, \mE\, T(\ux)  \hspace{8mm} \mbox{and} \hspace{8mm} \ux\, \uS_2 = -\, \Gamma T(\ux)
$$
this choice is not completely arbitrary since $\uS_1$ and $\uS_2$ always must satisfy the relation
\begin{equation}
\label{entangle}
\uS_1 + \uS_2 = \dirac_{rad}\, T(\ux) + \dirac_{ang}\, T(\ux) = \dirac\, T(\ux)
\end{equation}
where the distribution at the right--hand side is, quite naturally, uniquely determined once the distribution $T$ has been given.
One could say that the differential operators $\dirac_{rad}$ and $\dirac_{ang}$ are {\em entangled} in the sense that the results of their actions on a distribution $T$ are subject to (\ref{entangle}).\\
Let us give an simple example to illustrate this phenomenon. Consider the regular distribution $T(\ux) = \ux$. Then $\dirac\, \ux = -m$, $\mE\, \ux = \ux$ and $\Gamma\, \ux = (m-1)\, \ux$, whence
$$
\uom\, \p_r\, \ux = -1 + c_1\, \delta(\ux) \hspace{8mm} \mbox{and} \hspace{8mm} \frac{1}{r}\, \p_{\uom}\, \ux = 1 - m + c_2\, \delta(\ux)
$$
with the restriction that the arbitrary constants $c_1$ and $c_2$ always must satisfy the entanglement condition $c_1 + c_2 = 0$.\\

Apparently there seems to be no possibility to uniquely define the actions of the $\dirac_{rad}$ and $\dirac_{ang}$ operators on a standard distribution by singling out specific distributions in the equivalent classes (\ref{equirad}) and (\ref{equiang}), except for the following two special cases.\\

\noindent
(i) If the distribution $T(\ux)$ is {\em radial}, i.e. only depends on $r = |\ux|$: $T(\ux) = T^{rad}(|\ux|)$, then we put
$$
\frac{1}{r}\, \p_{\uom}\, T^{rad} = 0 \quad {\rm and} \quad  \uom\, \p_r\, T^{rad} = \dirac\, T^{rad}
$$ 
This first special case is illustrated by the delta-distribution (see also \cite{bsv}):
$
\frac{1}{r}\, \p_{\uom}\, \delta(\ux) = 0 $ and $\uom\, \p_r\, \delta(\ux) = \dirac\, \delta(\ux)
$.\\

\noindent
(ii) If the distribution $T(\ux)$ is {\em angular}, i.e. only depends on $\uom = \frac{\uom}{|\ux|}$, then we put
$$
\uom\, \p_r\, T = 0  \quad {\rm and} \quad  \frac{1}{r}\, \p_{\uom}\, T = \dirac\, T
$$
This second special case is illustrated by the regular distribution $\uom$ for which
$
\uom\, \p_r\, \uom = 0 $ and $ \frac{1}{r}\, \p_{\uom}\, \uom = \dirac\, \uom = -(m-1) \frac{1}{r}
$.\\

\noindent
In Section \ref{derivatives} we will point out two other cases where the actions of the $\dirac_{rad}$ and $\dirac_{ang}$ operators are uniquely defined.


\section{The Laplace operator in spherical co-ordinates}
\label{laplace}


As was already mentioned in the preceding section, the Dirac operator factorizes the Laplace operator: $-\, \Delta = \dirac^2$. As the Laplace operator is a scalar operator it holds that
$$
\Delta = - \, \dirac \cdot \dirac = |\dirac|^2
$$
Passing to spherical co-ordinates we obtain
\begin{align}
\Delta &= -\, (\dirac_{rad} + \dirac_{ang})^2 \nonumber \\
&= \p_r^2 + (m-1) \frac{1}{r}\, \p_r + \frac{1}{r^2}\, (\uom\, \p_{\uom} - \p_{\uom}^2) \nonumber\\
&= \p_r^2 + (m-1) \frac{1}{r}\, \p_r + \frac{1}{r^2}\, \Delta^* \nonumber
\end{align}
since
\begin{align}
\dirac_{rad} \, \dirac_{rad} &= -\, \p_{r}^2 \nonumber \\
\dirac_{rad} \, \dirac_{ang} &= -\frac{1}{r^2}\, \uom\, \p_{\uom} + \frac{1}{r}\, \uom\, \p_{\uom}\, \p_r \nonumber \\
\dirac_{ang} \, \dirac_{rad} &= -\, (m-1) \frac{1}{r}\, \p_r - \frac{1}{r}\, \p_r\, \uom\, \p_{\uom} \nonumber \\
\dirac_{ang} \, \dirac_{ang} &= \frac{1}{r^2}\, \p_{\uom}^2 \nonumber
\end{align}
The Laplace--Beltrami operator $\Delta^* = \uom\, \p_{\uom} - \p_{\uom}^2$, sometimes denoted by $\Delta_0$ and sometimes called the spherical Laplace operator, is a scalar angular operator; as $\uom\, \p_{\uom} = -\, \Gamma$ is a bivector operator, it follows that
$$
\Delta^* = -\, \p_{\uom} \cdot \p_{\uom} = |\p_{\uom}|^2 \qquad {\rm and} \qquad \uom\, \p_{\uom}  = \p_{\uom} \wedge \p_{\uom} = -\, \Gamma
$$
It is a nice observation that while the Laplace operator $\Delta$ is the normsquared of the Dirac operator, the Laplace--Beltrami operator is the normsquared of the spherical Dirac operator.\\

\noindent
As is the case for the Laplace operator $\Delta = \sum_{j=1}^m \, \p_{x_j}^2$, also the Lapalace--Beltrami operator may be expressed in terms of derivatives with respect to the cartesian co-ordinates.

\begin{proposition}
The angular differential operators $\p_{\uom}^2$ and $\Delta^*$ may be written in terms of cartesian derivatives as
$$
\p_{\uom}^2 = \Gamma^2 - (m-1)\, \Gamma
$$
and
$$
\Delta^* = (m-2)\, \Gamma - \Gamma^2
$$
\end{proposition}

\pf
One has
\begin{align}
\Gamma^2 &= (-\, \uom\, \p_{\uom})^2 =  \uom\, \p_{\uom}\, \uom\, \p_{\uom} \nonumber \\
&= \uom\, ( (1-m)   -    \uom\, \p_{\uom})  \p_{\uom} \nonumber \\
&= (1-m)\, \uom\, \p_{\uom} + \p_{\uom}^2 \nonumber \\
&= (m-1)\, \Gamma + \p_{\uom}^2 \nonumber
\end{align}
and
\begin{align}
\Delta^* &=  \uom\, \p_{\uom} - \p_{\uom}^2 \nonumber\\
&= -\, \Gamma - \Gamma^2 + (m-1)\, \Gamma \nonumber \\
&= (m-2)\, \Gamma -  \Gamma^2 \nonumber
\end{align} \eop

There is a second, and quite naturally equivalent, way to write the Laplace--Beltrami operator by means of cartesian derivatives. It only needs a straightforward calculation to prove the following result.

\begin{proposition}
One has
$$
\Delta^* = \sum_{j<k}\, L_{jk}^2 = \sum_{j<k}\, (x_j\, \p_{x_k} - x_k\, \p_{x_j})^2
$$
\end{proposition}

The actions of the Laplace operator and the Laplace--Beltrami operator on a distribution being uniquely well-defined, the question arises how to define the actions on a distribution of the three parts of the Laplace operator expressed in spherical co-ordinates. It turns out that these actions are well--defined, though not uniquely, through equivalent classes of distributions.

\begin{proposition}
\label{laplaceparts}
Let $T$ be a scalar distribution. One has
\begin{itemize}
\item[(i)] $\p_r^2\, T = S_2 + \delta(\ux)\, c_2 - \sum_{j=1}^m\, c_{1,j}\, \p_{x_j} \delta(\ux)$\\[1mm] for arbitrary constants $c_2$ and $ c_{1,j}, j=1,\ldots,m$ and any distribution $S_2$ such that $\ux\, S_2 = \mE\, \underline{S}_1$ with $\ux\, \underline{S}_1 = -\, \mE\, T$
\item[(ii)] $\invr\, \p_r\, T = S_3 + \frac{1}{m}\, \sum_{j=1}^m\, c_{1,j}\, \p_{x_j} \delta(\ux) + c_3\, \delta(\ux)$\\[1mm] for arbitrarily constant $c_3$ and any distribution $S_3$ such that $\ux\, S_3 = \underline{S}_1$
\item[(iii)] $\frac{1}{r^2}\ \Delta^*\, T = S_4 + c_4\, \delta(\ux) + \sum_{j=1}^m\, c_{5,j}\, \p_{x_j} \delta(\ux)$\\[1mm] for arbitrary constants $c_4$ and $ c_{5,j}, j=1,\ldots,m$ and any distribution $S_4$ such that $r^2\, S_4 = \Delta^*\, T$
\end{itemize}
\end{proposition}

\pf
\begin{itemize}
\item[(i)] From Section \ref{diracoperator} we know that
$$
(\uom\, \p_r)\, T = -\, \left[ \frac{1}{\ux}\, \mE\, T  \right] = \underline{S}_1 + \delta(\ux)\, \underline{c}_1 
$$
with $\ux\, \underline{S}_1 = -\, \mE\, T$. It follows that
\begin{align}
\p_r^2\, T &= -\, (\uom\, \p_r)^2\, T   \nonumber\\
&= -\, (\uom\, \p_r)\, (\underline{S}_1 + \delta(\ux)\, \underline{c}_1) \nonumber\\
&= \left[ \frac{1}{\ux}\, \mE\, \underline{S}_1  \right] - \dirac \delta(\ux)\, \underline{c}_1 \nonumber\\
&= S_2 + \delta(\ux)\, c_2 - \dirac \delta(\ux)\, \underline{c}_1 \nonumber
\end{align}
with $\ux\, S_2 = \mE\, \underline{S}_1$.
\item[(ii)] We have consecutively
\begin{align}
\invr\, \p_r\, T &= \frac{1}{\ux}\, (\uom\, \p_r)\, T \nonumber\\
&= \frac{1}{\ux}\, ( \underline{S}_1 + \delta(\ux)\, \underline{c}_1 ) \nonumber\\
&= S_3 +  \frac{1}{\ux}\, \delta(\ux)\, \underline{c}_1 \nonumber\\
&= S_3 + \frac{1}{m}\, \dirac\, \delta(\ux)\,  \underline{c}_1+ \delta(\ux)\, c_3 \nonumber
\end{align}
with $\ux\, S_3 = \underline{S}_1$.
\item[(iii)] The distribution $\Delta^*\, T$ is uniquely defined and $r^2 $ is an analytic function with a second order zero at the origin. The result follows immediately.
\end{itemize}

\begin{remark}
{\rm
The operators $\p_r^2$,  $\invr\, \p_r$ and $\frac{1}{r^2}\ \Delta^*$ are {\em entangled} in the sense that, given a distribution $T$ and having chosen appropriately the distributions $\underline{S}_1$, $S_2$, $S_3$ and $S_4$, all arbitrary constants appearing in the expressions of Proposition \ref{laplaceparts} should satisfy the entanglement condition  
$$
\p_r^2\, T + (m-1) \frac{1}{r}\, \p_r\, T + \frac{1}{r^2}\, \Delta^*\, T = \Delta\, T
$$ 
the distribution at the right--hand side being uniquely determined.
}
\end{remark}

\begin{example}
{\rm
Proposition \ref{laplaceparts} may be generalised to distributions which are e.g. vector valued. Let us illustrate this by considering the distribution $T = \ux^3 = -\, r^3\, \uom$, for which, by a direct computation, $\Delta\, T = \Delta(\ux^3) = -\, 2\, (m+2)\, \ux$, and $\Delta^*\, T = \Delta^*\, (\ux^3) = (m-1)\, r^2\, \ux = -\, (m-1)\, \ux^3$.\\
As $\mE\, T = \mE\, (\ux^3) = 3\, \ux^3$, we chose $S_1 = -\, 3\, \ux^2 = 3\, r^2$ satisfying $\ux\, S_1 = -\, 3\, \ux^3$. As $\mE\, S_1 = \mE\, (-3\ux^2) = -\, 6\, \ux^2 = 6\, r^2$, we chose $\underline{S}_2 = -\, 6\, \ux$ satisfying $\ux\, \underline{S}_2 = -\, 6\, \ux^2$, and $\underline{S}_3 = -\, 3\, \ux$ satisfying $\ux\, \underline{S}_3 = -\, 3\, \ux^2$. Finally we chose $\underline{S}_4 = (m-1)\, \ux$, satisfying $r^2\, \underline{S}_4 = \Delta^*\, T = (m-1)\, r^2\, \ux$.\\
This leads to:
\begin{itemize}
\item[(i)] $\p_r^2\, T = \p_r^2\, (\ux^3) = -\, 6\, \ux + \delta(\ux)\, c_2 - \sum_{j=1}^m\, c_{1,j}\, \p_{x_j} \delta(\ux)$

\item[(ii)] $\invr\, \p_r\, T = \invr\, \p_r\, (\ux^3) = -\, 3\, \ux + \frac{1}{m}\, \sum_{j=1}^m\, c_{1,j}\, \p_{x_j} \delta(\ux) + c_3\, \delta(\ux)$

\item[(iii)] $\frac{1}{r^2}\, \Delta^*\, T = \frac{1}{r^2}\, \Delta^*\, (\ux^3) = (m-1)\, \ux + c_4\, \delta(\ux) + \sum_{j=1}^m\, c_{5,j}\, \p_{x_j} \delta(\ux)$

\end{itemize}
provided that the arbitrary constants should satisfy the entanglement conditions
$$
\begin{cases}
c_2 + (m-1)\, c_3 + c_4 = 0\\
-\, \frac{1}{m}\, c_{1,j} + c_{5,j} = 0, \ \  j=1,\ldots,m
\end{cases}
$$
}
\end{example}

\begin{remark}
{\rm
In the special case where the distribution $T^{rad}$ is radial: $T^{rad}(\ux) = T^{rad}(|\ux|)$, it holds that $\Delta^*\, T^{rad} = 0$ and
$$
\Delta\, T^{rad} = \p^2_{r}\, T^{rad} + (m-1)\, \invr\, \p_r\, T^{rad}
$$
Invoking the SO$(m)$--invariance of the distributions $T^{rad}$, $\p^2_{r}\, T^{rad}$ and $\invr\, \p_r\, T^{rad}$, it then holds that 
$$
\p^2_{r}\, T^{rad} = S_2 + c_2\, \delta(\ux), \quad \ux\, S_2 = \mE\, \underline{S}_1, \quad \ux\, \underline{S}_1 = -\, \mE\, T^{rad}
$$
and
$$
\invr\, \p_r\, T^{rad} = S_3 + c_3\, \delta(\ux), \quad \ux\, S_3 = \underline{S}_1
$$
where the arbitrary constants $c_2$ and $c_3$ must satisfy the entanglement condition
$$
S_2 + (m-1)\, S_3 + (c_2 + (m-1)\, c_3)\, \delta(\ux) = \Delta\, T^{rad}
$$
the distributions $S_2$ and $S_3$ having been chosen appropriately.\\
As an example consider the radial distribution $T^{rad} = r$ for which $\Delta\, r = (m-1)\invr$. Then $ \p^2_{r}\, r = c_2\, \delta(\ux)$ and $ \invr\, \p_r\, r = \invr + c_3\, \delta(\ux)$ where the arbitrary constants must satisfy the entanglement condition $c_2 + (m-1)\, c_3 = 0$.\\
Finally notice that in the very special case of the delta-distribution $\delta(\ux)$ it was proved in \cite{bsv}, invoking the homogeneity of $\delta(\ux)$, that 
$$
\p^2_{r}\, \delta(\ux) = \onehalf\, (m+1)\, \Delta\, \delta(\ux)
$$
and
$$
\invr\, \p_r\, \delta(\ux) = -\, \onehalf\, \Delta\, \delta(\ux)
$$
}
\end{remark}


\section{Radial and angular derivatives of distributions}
\label{derivatives}


In Section 1 we explained why it is impossible to define within the class of distributions the radial derivative $\p_{r}\, T$ and the vector angular derivative $\p_{\uom}\, T$ of a distribution $T$. Neither is it possible to multiply a distribution by the non--analytic functions $r$ and $\uom$. For {\em legitimizing} those {\em forbidden actions} we have to take the signumdistributions into consideration instead. 

\begin{definition}
\label{multiplicomega}
The product of a scalar--valued distribution $T$ by the function $\uom$ is the signumdistribution $T^{\vee}$ associated to $T$, and it holds for all test functions $\uom\, \varphi \in \Omega(\mR^m, \mR^m)$ that
$$
\langle \ \uom\, T \ , \ \uom\, \varphi \ \rangle = \langle \ T^{\vee} \ , \ \uom\, \varphi \ \rangle = -\, \langle \  T \  ,  \varphi \ \rangle
$$
\end{definition}

\begin{example}
{\rm
The product of the delta-distribution $\delta(\ux)$ by the function $\uom$ is the signumdistribution $\delta(\ux)^{\vee} = \uom\, \delta(\ux)$ for which
$$
\langle \, \uom\, \delta(\ux) \, , \, \uom\, \varphi \, \rangle = -\, \langle \, \delta(\ux) \, , \,  \varphi \, \rangle = -\, \varphi(0)
$$
for all test functions $\uom\, \varphi \in \Omega(\mR^m, \mR^m)$.
}
\end{example}

\begin{definition}
\label{multiplicr}
The product of a scalar--valued distribution $T$ by the function $r$ is the signumdistribution $r\,T = (-\, \ux\ T)^{\vee}$ given by
$$
\langle \ r\, T \ , \ \uom\, \varphi \ \rangle =  \langle \  \ux\, T \  , \ \varphi \ \rangle =  \langle \  T \  , \  \ux\, \varphi \ \rangle
$$
according to (the boldface part of) the commutative scheme
$$
\begin{array}{ccccccccc}
  &&    & {}_{-\ux}   &                             &  &&&\\
  && {\bf T} & \longrightarrow  & {\bf -\, \ux\, T} & &&& \\[1mm]
  &&\hspace{-8mm}{}_{- \uom} &\hspace{2mm}{} \hspace{11mm}{}_{- r}&\hspace{6mm}{}_{-\uom}&&&&\\[-2mm]
  && \uparrow & \diagdown \hspace{-1.3mm} \nearrow & \uparrow &&&&\\[-1.1mm]
  && \downarrow & \diagup \hspace{-1.3mm} \searrow & \downarrow &&&&\\[-2mm]
  &&\hspace{-4.5mm}{}^{\uom}&\hspace{2mm}{} \hspace{10mm}{}^{r}&\hspace{5mm}{}^{\uom}&&&&\\[1mm]
  && T^{\vee} = \uom\, T & \longrightarrow & {\bf r\, T} &&&&\\[-1mm]
  &&               & {}_{-\ux}  &           &&&&
 \end{array}
$$
\end{definition}

\begin{example}
{\rm
The product of the delta-distribution $\delta(\ux)$ by the function $r$ is the zero distribution since $(-\ux)\, \delta(\ux) = 0$. The product of the distribution $\dirac\, \delta(\ux)$ by $r$ is the signumdistribution $(-m)\, \uom\, \delta(\ux)$ since $(-\ux)\, \dirac\, \delta(\ux) = \mE\, \delta(\ux) = (-m)\, \delta(\ux)$.
}
\end{example}

\begin{remark}
{\rm
In the commutative scheme of Definition \ref{multiplicr}, and in all the commutative schemes in the sequel of this paper as well, the upper row is situated in the world of distributions, while the objects in the lower row are signumdistributions. Vertical transition from the distributions to the signumdistributions and vice versa is executed by the multiplication operators $\uom$ and $-\, \uom$ respectively. Each of the horizontally acting operators between distributions, has its counterpart in the world of signumdistributions, and vice versa; e.g. in the above commutative scheme the multiplication operator $-\, \ux$ between the distributions $T$ and $-\, \ux\, T$ corresponds with the multiplication operator $-\, \ux$ between the signumdistributions $T^{\vee}$ and $(-\, \ux\, T)^{\vee} = r\, T$. In fact this implies the definition of the multiplication of the signumdistribution $T^{\vee} = \uom\, T$ by the function $\ux$ resulting in the signumdistribution $-\, r\, T$. This idea will be taken up again in the next section.
}
\end{remark}

\begin{definition}
\label{radderiv}
The derivative with respect to the radial distance $r$ of a scalar--valued distribution $T$ is the equivalent class of signumdistributions 
$$
\left[ \p_r\,T \right] = \left[ -\, \uom\, \p_r\, T \right]^{\vee} = \left[\frac{1}{\ux} \, \mE\, T \right]^{\vee} = \left( S + \underline{c}\, \delta(\ux) \right)^{\vee} =  \uom\, S + \uom\, \delta(\ux)\, \underline{c}
$$ 
for any vector distribution $S$ satisfying $\ux \, S =  \mE\, T$, according to (the boldface part of) the commutative scheme
$$
\begin{array}{ccccccccc}
  &&    & {}_{-\uom \p_r}   &                             &&&&\\
  && {\bf T} & \longrightarrow  & {\bf  \left[\frac{1}{\ux} \, \mE\, T \right] }&&&& \\[1mm]
  &&\hspace{-8mm}{}_{- \uom} &\hspace{2mm}{} \hspace{11mm}{}_{-\p_r}&\hspace{6mm}{}_{-\uom}&&&&\\[-2mm]
  && \uparrow & \diagdown \hspace{-1.3mm} \nearrow & \uparrow &&&&\\[-1.1mm]
  && \downarrow & \diagup \hspace{-1.3mm} \searrow & \downarrow &&&&\\[-2mm]
  &&\hspace{-4.5mm}{}^{\uom}&\hspace{2mm}{} \hspace{10mm}{}^{\p_r}&\hspace{5mm}{}^{\uom}&&&&\\[1mm]
  && T^{\vee} = \uom\, T & \longrightarrow & {\bf \left[ \p_r T \right] } &&&&\\[-1mm]
  &&               & {}_{-\uom \p_r}  &           &&&&
 \end{array}
$$
\end{definition}

\begin{remark}
{\rm
In the special case of a scalar--valued {\rm radial} distribution $T^{rad}$, its radial derivative $\p_r\,T^{rad}$ is uniquely determined to be the signumdistribution $\p_r\,T^{rad} = (-\, \dirac\, T^{rad})^{\vee}$ given for all test functions $\uom\, \varphi$ by
$$
\langle \ \p_r\, T^{rad} \ , \  \uom\, \varphi \ \rangle =  \langle \  \uom \p_r \, T^{rad} \  , \  \varphi \ \rangle = \langle \  \dirac \, T^{rad} \  ,  \ \varphi \ \rangle
$$
according to (the boldface part of) the commutative scheme
$$
\begin{array}{ccccccccc}
  &&    & {}_{-\uom \p_r}   &                             &&&&\\
  && {\bf T^{rad}} & \longrightarrow  & {\bf  -\dirac\, T^{rad}}&&&& \\[1mm]
  &&\hspace{-8mm}{}_{- \uom} &\hspace{2mm}{} \hspace{11mm}{}_{-\p_r}&\hspace{6mm}{}_{-\uom}&&&&\\[-2mm]
  && \uparrow & \diagdown \hspace{-1.3mm} \nearrow & \uparrow &&&&\\[-1.1mm]
  && \downarrow & \diagup \hspace{-1.3mm} \searrow & \downarrow &&&&\\[-2mm]
  &&\hspace{-4.5mm}{}^{\uom}&\hspace{2mm}{} \hspace{10mm}{}^{\p_r}&\hspace{5mm}{}^{\uom}&&&&\\[1mm]
  && \uom\, T^{rad} & \longrightarrow & {\bf \p_r T^{rad} = -\uom\, \dirac\, T^{rad} } &&&&\\[-1mm]
  &&               & {}_{-\uom \p_r}  &           &&&&
 \end{array}
$$
}
\end{remark}

\begin{example}
{\rm
The radial derivative of the delta-distribution is the signumdistributioin $\p_r\, \delta(\ux) = (-\, \dirac\, \delta(\ux))^{\vee} = -\, \uom\, \dirac\, \delta(\ux)$.
}
\end{example}

\begin{remark}
{\rm
The commutative scheme of Definition \ref{radderiv} implies the definition of the action of the operator $\dirac_{rad} = \uom\, \p_r$ on the signumdistribution $T^{\vee} = \uom\, T$ resulting in the (equivalence class of) signumdistributions $-\, [\p_r\, T]$.
In the special case where the distribution $T$ is radial: $T = T^{rad}$, the action result of the operator $\dirac_{rad} = \uom\, \p_r$ on $\uom\, T^{rad}$ is the uniquely determined signumdistribution
$$
(\uom\, \p_r) \uom\, T^{rad} = -\, \p_r\, T^{rad} = \uom(\uom\, \p_r)T = \uom\, \dirac\, T = (\dirac\, T)^{\vee}
$$
and it then holds for all test functions $\uom\, \varphi$ that
$$
\langle \  -\, \uom\, \p_r\, T^{\vee} \ , \ \uom\, \varphi \ \rangle = \langle \  ( \uom\, \p_r\, T^{\vee})^{\wedge}  \ , \   \varphi \ \rangle = \langle \ \p_r\, T^{\vee} \ , \ \varphi \ \rangle =
\langle \ -\, (\p_r\, T)^{\wedge} \ , \ \varphi \ \rangle = \langle \ \p_r\, T \ , \ \uom\, \varphi \ \rangle
$$ 
}
\end {remark}

\begin{definition}
\label{angderiv}
The angular $\p_{\uom}$--derivative of a scalar--valued distribution $T$ is the signumdistribution $\p_{\uom}\,T = (\Gamma\, T)^{\vee}$ given for all test functions $\uom\, \varphi$ by
$$
\langle \  \uom\, \varphi\ , \  \p_{\uom}\, T \ \rangle =  \langle \   \varphi  \  , \  \uom\, \p_{\uom} \, T  \ \rangle = \langle \  \varphi  \  ,  \  -\, \Gamma \, T \ \rangle
$$
according to (the boldface part of) the commutative scheme
$$
\begin{array}{ccccccccc}
  &&    & {}_{- \uom \p_{\uom}}   &                             &  &&&\\
  && {\bf T} & \longrightarrow  & {\bf \Gamma\, T} & &&& \\[1mm]
  &&\hspace{-8mm}{}_{- \uom} &\hspace{2mm}{} \hspace{13mm}{}_{\uom \p_{\uom} \uom}&\hspace{6mm}{}_{-\uom}&&&&\\[-2mm]
  && \uparrow & \diagdown \hspace{-1.3mm} \nearrow & \uparrow &&&&\\[-1.1mm]
  && \downarrow & \diagup \hspace{-1.3mm} \searrow & \downarrow &&&&\\[-2mm]
  &&\hspace{-4.5mm}{}^{\uom} & \hspace{2mm}{} \hspace{10mm}{}^{\p_{\uom}}&\hspace{5mm}{}^{\uom}&&&&\\[1mm]
  && T^{\vee} = \uom\, T & \longrightarrow & {\bf \p_{\uom}\, T}  &&&&\\[-1mm]
  &&               & {}_{- \p_{\uom} \uom}  &           &&&&
 \end{array}
$$
\end{definition}

\begin{example}
{\rm
For the delta-distribution it holds that $\p_{\uom}\, \delta(\ux) = 0$ since $\Gamma\, \delta(\ux) = 0$.
}
\end{example}

\begin{remark}
{\rm
The commutative scheme of Definition \ref{angderiv} implies the definition of the action of the operator $ \p_{\uom}\, \uom$ on the signumdistribution $T^{\vee} = \uom\, T$ resulting in the signumdistribution $-\, \p_{\uom}\, T$, which in its  turn implies the definition of the action of the $\Gamma$--operator on the signumdistribution $T^{\vee} = \uom\, T$ resulting in the signumdistribution
$$
\Gamma(\uom\, T) = (m-1)\, \uom\, T - \p_{\uom}\, T
$$
since
$
\p_{\uom}\, \uom = (1-m){\bf 1} - \uom\, \p_{\uom} = (1-m){\bf 1} + \Gamma
$.
}
\end{remark}

  
\section{Actions on signumdistributions}
\label{actionsignum}
  

Firstly let us make the following observations. With each well--defined operator $P$ acting between distributions there corresponds an operator $P^{\vee} = \uom\, P\, (-\uom)$ acting between signumdistributions according to the following commutative scheme
$$
\begin{array}{ccccccccc}
  &&    & {}_{P}   &                             &&&&\\
  &&  ^{s}U^{\wedge} = -\uom\,  ^{s}U & \longrightarrow  &  P\,   (-\uom)\,  ^{s}U &&&& \\[1mm]
  &&\hspace{-8mm}{}_{- \uom}&&&&&&\\[-2mm]
  && \uparrow &   & | &&&&\\[-1.1mm]
  && | &   & \downarrow &&&&\\[-2mm]
  &&&&\hspace{5mm}{}^{\uom}&&&&\\[1mm]
  && ^{s}U & \longrightarrow & P^{\vee}\, ^{s}U = \uom\, P\, (-\uom)\, ^{s}U &&&&\\[-1mm]
  &&               & {}_{P^{\vee}}  &           &&&&
 \end{array}
$$
in this way giving rise to a pair of operators $P$ and $P^{\vee}$ which we call a {\em signum--pair of operators}. If the action result of the operator $P$ is uniquely determined then the action result of $P^{\vee}$ is uniquely determined too, in which case we use the notation $(P,P^{\vee})$ for this signum--pair of operators. If, on the contrary, the action result of $P$ is an equivalence class of distributions, then the action result of $P^{\vee}$ will be an equivalence class of signumdistributions as well, in which case we use the notation $[P, P^{\vee}]$. Whenever one has a signum--pair of operators $(P,P^{\vee})$ or $[P, P^{\vee}]$ at one's disposal, one defines the action of $P^{\vee}$ on a signumdistribution as follows.

\begin {definition}
\label{actionsignumdistrib}
If the operator $P^{\vee}$ is the signum--partner to the operator $P$ acting between distributions, then the action of  $P^{\vee}$ on a signumdistribution $^{s}U$ is given by
$$
P^{\vee}\, ^{s}U = \left(P\, (^{s}U)^{\wedge}  \right)^{\vee} = \uom\, P\, (-\uom)^{s}U
$$
\end{definition}

\begin{remark}
{\rm
When the operators $P$ and $P^{\vee}$ form a signum--pair of operators then it follows that $P = (-\uom)\, P^{\vee}\, \uom$, and it becomes tempting to consider the operator $P$ as the signum--partner to the operator $P^{\vee}$, and to define the action of the operator $P$ on signumdistributions through the action of $P^{\vee}$ on distributions. However this is legitimate only if the operator $P^{\vee}$ is a ''legal`` well--defined operator between distributions.
}
\end{remark}

\noindent
The following pairs of operators, some of which we already encountered in the previous sections, may easily be checked to be signum--pairs: $(\ux, \ux)$,  $(r^2, r^2)$, $[\uom\, \p_r, \uom\, \p_r]$ and $(\mE, \mE)$, which implies the definition of the action of the operators $\ux$, $r^2$, $\uom\, \p_r$ and $\mE$  on a signumdistribution following Definition \ref{actionsignumdistrib}.\\
We also found (see Definition \ref{angderiv}) the signum--pair $(\Gamma \, , -\, \p_{\uom}\, \uom) = (\Gamma,  (m-1){\bf 1} -\, \Gamma)$. As the operator $(m-1){\bf 1} -\, \Gamma$ clearly is a well--defined operator acting on distributions, also $$( -\, \p_{\uom}\, \uom \, , \, \Gamma) = ((m-1){\bf 1} -\, \Gamma \, , \, \Gamma)$$ is a signum--pair of operators, which implies the definition of the action of the $\Gamma$ operator on signumdistributions through Definition \ref{actionsignumdistrib}. A straightforward computation then yields the additional signum--pairs of operators $(\Gamma^2  \, , \, \Gamma^2 - 2(m-1)\, \Gamma + (m-1)^2{\bf 1})$ and its symmetric\\  $(\Gamma^2 - 2(m-1)\, \Gamma + (m-1)^2{\bf 1} \, , \, \Gamma^2)$.\\

From Section \ref{laplace}, Proposition 2 we know the differential operators $\p_{\uom}^2$ and $\Delta^*$ to be ''cartesian``, i.e. they may be expressed in terms of derivatives with respect to the cartesian co-ordinates. Let us compute their signum--partners. We find
\begin{eqnarray*}
\uom\, \p_{\uom}^2\, (-\, \uom) &=& \uom\, \p_{\uom}\, ( \uom\, \p_{\uom} + (m-1))\\
&=& \uom\, (-\ \uom\, \p_{\uom} - (m-1))\, \p_{\uom}  + (m-1)\, \uom\, \p_{\uom}\\
&=& \p_{\uom}^2
\end{eqnarray*}
and
\begin{eqnarray*}
\Delta^{* \vee} &=& (\uom\, \p_{\uom} - \p_{\uom}^2)^{\vee}\\
&=& - \p_{\uom}^2 - \uom\, \p_{\uom} - (m-1)\\
&=& -\, \Gamma^2 + m\, \Gamma - (m-1)\,{\bf 1}
\end{eqnarray*}
leading to the signum--pairs of operators $(\p_{\uom}^2 \, , \, \p_{\uom}^2)$, $( \Delta^* \, , \, {\bf Z}^* )$ and $( {\bf Z}^* \, , \, \Delta^* )$, where we have introduced the notation ${\bf Z}^* = \Delta^{* \vee}$. These signum--pairs of operators imply the definition of the actions of the operators $\p_{\uom}^2$, ${\bf Z}^*$ and $\Delta^{*}$ on signumdistributions.\\

Now let us determine the signum--partner $\underline{D}$ of the Dirac operator $\dirac$. We obtain
\begin{eqnarray*}
\underline{D} &=& \uom\, \dirac\, (-\uom) =  \uom\, (\uom\, \p_r + \invr\, \p_{\uom})\, (-\uom)\\
                      &=& \uom\, \p_r + \invr\, \uom\, \p_{\uom}\, (-\uom)\\
                      &=& \uom\, \p_r - \invr\,  \p_{\uom} + (m-1)\, \invr\, \uom
\end{eqnarray*}
giving rise to the signum--pairs of operators $(\dirac\, ,\, \underline{D})$ and $[\invr\, \p_{\uom} \, , \, - \invr\,  \p_{\uom} + (m-1)\, \invr\, \uom]$. The actions of the operators $\uD$ and $- \invr\,  \p_{\uom} + (m-1)\, \invr\, \uom$ on a signumdistribution then are defined along the lines of Definition \ref{actionsignumdistrib}.\\
Notice that while the actions of the operator $\dirac$ on distributions and of its signum--partner $\underline{D}$ on signumdistributions are uniquely defined, the action results of the operator $\invr\, \p_{\uom}$ on distributions and of its signum--partner $- \invr\, \p_{\uom} + (m-1)\, \invr\, \uom$ on signumdistributions are equivalent classes, both these operators remaining entangled with the operator $\uom\, \p_r$. Let us illustrate these observations by the following simple example. While $\dirac\, (-\ux) = m$ and $\underline{D}\, (-x)^{\vee} = \underline{D}\, r = m\, \uom$, it holds that $\uom\, \p_r\, (-\ux) = [\, 1\, ], \, \uom\, \p_r\, r = [\, \uom\, ]$ and $\invr\, \p_{\uom}\, (-\ux) = [\, m-1\, ], \left(\invr\, \p_{\uom} \right)^{\vee}\, r = [\, (m-1)\, \uom\, ]$, where the arbitrary constants appearing in these equivalence classes have to be chosen in order to fulfil the entanglement conditions $[\, 1\, ] + [\, m-1\, ] = m$ and $[\, \uom\, ] + [\, (m-1)\, \uom\, ] = m\, \uom$ respectively.\\
We could call the operator $\uD$ the {\em signum--Dirac operator}. In the same way as the Dirac operator factorizes the Laplace operator: $\dirac^2 = -\, \Delta$, the signum--Dirac operator factorizes the {\em signum--Laplace operator}, i.e. the signum--partner of the Laplace operator:
$$
\uD^2 = \left(\uom\, \dirac\, (-\uom)\right)^2 = \uom\, \dirac^2 \, (-\, \uom) = -\, \uom\, \Delta\, (-\, \uom) = - \Delta^{\vee}
$$
Introducing the notation ${\bf Z} = \Delta^{\vee}$, it follows that $(\Delta\, , \, {\bf Z})$ is a signum--pair of operators, with
\begin{eqnarray*}
{\bf Z} &=& -\, \uD^2\\
&=& -\, \left( \uom\, \p_r - \invr\,  \p_{\uom} + (m-1)\, \invr\, \uom  \right)^2\\
&=& \p_r^2 + (m-1)\, \invr\, \p_r - \invrsq\, ( \p_{\uom}^2 + \uom\, \p_{\uom} + (m-1))\\
&=& \p_r^2 + (m-1)\, \invr\, \p_r + \invrsq\, {\bf Z}^*
\end{eqnarray*}
It is not clear if the actions of the Dirac operator $\dirac$ and the Laplace operator $\Delta$ on signumdistributions may be defined through the possible signum--pairs of operators $(\uD \, , \, \dirac)$ and $({\bf Z} \, , \, \Delta)$, since for the moment the operators $\uD$ and ${\bf Z}$ are known only in terms of spherical co-ordinates. We postpone this discussion to Section \ref{division}.\\

Similarly with each well--defined operator $Q$ mapping a distribution to a signumdistribution, there corresponds a well--defined operator $Q^c = (-\uom)\, Q\, (-\uom) =  \uom\, Q\, \uom$ mapping a signumdistribution to a distribution according to the following commutative scheme
$$
\begin{array}{ccccccccc}
  &&    & {}_{P}   &                             &&&&\\
  &&  ^{s}U^{\wedge} = -\uom\,  ^{s}U & \longrightarrow  &  Q^c \, ^{s}U &&&& \\[1mm]
  &&\hspace{-8mm}{}_{- \uom} &\hspace{2mm}{} \hspace{11mm}{}_{Q^c}&\hspace{6mm}{}_{-\uom}&&&&\\[-2mm]
  && \uparrow & \diagdown \hspace{-1.3mm} \nearrow & \uparrow &&&&\\[-1.1mm]
  && | & \diagup \hspace{-1.3mm} \searrow & | &&&&\\[-2mm]
  &&&\hspace{2mm}{} \hspace{12mm}{}^{Q}&&&&&\\[1mm]
  &&   ^{s}U & \longrightarrow & Q \, (-\uom)\, ^{s}U &&&&\\[-1mm]
  &&               & {}_{P^{\vee}}  &           &&&&
 \end{array}
$$
in this way giving rise to a pair of operators $Q$ and  $Q^c$ which we call a {\em cross--pair of operators} denoted by either $(Q\, , Q^c)$ or $[Q\, , Q^c]$ depending on the nature of their action result similarly as in the case of a signum--pair of operators. That action of the operator $Q^c$ mapping a signumdistribution to a distribution then is defined as follows.

\begin{definition}
\label{actioncross}
Let $Q$ be a well--defined operator mapping a distribution to a signumdistribution, then the action of its cross--partner $Q^c$ on a signumdisribution $^{s}U$ is defined by $Q^c\, ^{s}U = Q\, (^{s}U^{\wedge}) = Q\, (-\uom)\, ^{s}U$.
\end{definition}

\noindent
In the previous sections we already encountered several cross--pairs of operators such as $(\uom\, , -\uom)$, $(r\, , -r)$, $[\p_r\, , -\p_r]$ and $(\p_{\uom}\, , \uom\, \p_{\uom}\, \uom)$, which in fact imply the definition of the action on signumdistributions of the operators $-\, \uom$, $-\, r$
, $-\, \p_r$ and  $\uom\, \p_{\uom}\, \uom$ respectively along the lines of Definition \ref{actioncross}. Let us state these definitions explicitly.

\begin{definition}
\label{signummultiplicomega}
The product of a scalar--valued signumdistribution $^{s}U$ by the function $\uom$ is the distribution $-\, ^{s}U^{\wedge}$ associated to the signumdistribution $-\, ^{s}U$, and it holds for all test functions $\varphi$ that
$$
\langle \ \uom\, ^{s}U \ , \  \varphi \ \rangle = \langle \ -\, ^{s}U^{\wedge} \ , \  \varphi \ \rangle =  \langle \  ^{s}U \  ,  \ \uom\, \varphi \ \rangle
$$
\end{definition}

\begin{definition}
\label{signummultiplicr}
The product of a scalar--valued signumdistribution $^{s}U$ by the function $r$ is the distribution $r\, ^{s}U = \ux\, (^{s}U)^{\wedge}$ given for all test functions $\varphi$ by
$$
\langle \ r\, ^{s}U \ , \  \varphi \ \rangle =  \langle \   \ux\, (-\, \uom\, ^{s}U)  \  , \  \varphi \ \rangle =  \langle \  -\,  \uom\, ^{s}U \  , \  \ux\, \varphi \ \rangle = \langle \   ^{s}U \  , \  -\, \uom\, (\ux\, \varphi) \ \rangle
$$
according to (the boldface part of) the commutative scheme
$$
\begin{array}{ccccccccc}
  &&    & {}_{\ux}   &                             &&&&\\
  && {\bf ^{s}U^{\wedge} = -\uom\,  ^{s}U} & \longrightarrow  & {\bf r \, ^{s}U} &&&& \\[1mm]
  &&\hspace{-8mm}{}_{- \uom}&\hspace{2mm}{} \hspace{11mm}{}_{r}&\hspace{6mm}{}_{-\uom}&&&&\\[-2mm]
  && \uparrow & \diagdown \hspace{-1.3mm} \nearrow & \uparrow &&&&\\[-1.1mm]
  && \downarrow & \diagup \hspace{-1.3mm} \searrow & \downarrow &&&&\\[-2mm]
  &&\hspace{-4.5mm}{}^{\uom}&\hspace{2mm}{} \hspace{12mm}{}^{- r}&\hspace{5mm}{}^{\uom}&&&&\\[1mm]
  &&  {\bf ^{s}U} & \longrightarrow & \ux  \, ^{s}U &&&&\\[-1mm]
  &&               & {}_{\ux}  &           &&&&
 \end{array}
$$
\end{definition}

\begin{definition}
\label{signumradderiv}
The derivative with respect to the radial distance $r$ of a scalar--valued signumdistribution $^{s}U$ is the equivalent class of distributions
$$
\left[ \p_r\,^{s}U \right] = \left [  \uom\, \p_r\, ^{s}U^{\wedge}  \right] = \left[-\, \frac{1}{\ux} \, \mE\, ^{s}U^{\wedge}  \right] =  \left[\frac{1}{\ux} \, \mE\, \uom\, ^{s}U  \right] = S + c\, \delta(\ux) 
$$
for any scalar distribution $S$ satisfying $\ux \, S =  -\, \mE\, ^{s}U^{\wedge} = \mE\, \uom\, ^{s}U$, 
according to (the bold face part of) the commutative scheme
$$
\begin{array}{ccccccccc}
  &&    & {}_{\uom \p_r}   &                             &&&&\\
  && {\bf ^{s}U^{\wedge} = -\uom\,  ^{s}U} & \longrightarrow  & {\bf \left[ \p_r \, ^{s}U \right] }&&&& \\[1mm]
  &&\hspace{-8mm}{}_{- \uom}&\hspace{2mm}{} \hspace{11mm}{}_{\p_r}&\hspace{6mm}{}_{-\uom}&&&&\\[-2mm]
  && \uparrow & \diagdown \hspace{-1.3mm} \nearrow & \uparrow &&&&\\[-1.1mm]
  && \downarrow & \diagup \hspace{-1.3mm} \searrow & \downarrow &&&&\\[-2mm]
  &&\hspace{-4.5mm}{}^{\uom}&\hspace{2mm}{} \hspace{12mm}{}^{-\p_r}&\hspace{5mm}{}^{\uom}&&&&\\[1mm]
  && {\bf ^{s}U} & \longrightarrow & \uom\, \left[\p_r\, ^{s}U \right] &&&&\\[-1mm]
  &&               & {}_{\uom \p_r}  &           &&&&
 \end{array}
$$
\end{definition}

\begin{remark}
{\rm
As we have now at our disposal the definitions of the multiplication by $r$ (Definition \ref{signummultiplicr}) and of the radial derivative $\p_r$ (Definition \ref{signumradderiv}) of a signumdistribution, we can reconsider the action of the Euler operator $\mE = r\, \p_r$ on the signumdistribution $^{s}U$, resulting into the unique signumdistribution $\mE\, ^{s}U$ given by
$$
\mE\, ^{s}U = (r\, \p_r)\, ^{s}U = r\, (\p_r\, ^{s}U) = \uom\, (-\ux\, [\p_r\, ^{s}U]) = \uom\, (\mE\, ^{s}U^{\wedge}) = \uom\, (-\ux\, S) = r\, S
$$
for any distribution $S$ satisfying $\ux\, S = -\, \mE\, ^{s}U^{\wedge}$, according to the commutative scheme:
$$
\begin{array}{ccccccccc}
  &&    & {}_{\uom \p_r}   &                             &  -\ux &&\\
  &&  ^{s}U^{\wedge} = -\uom\,  ^{s}U & \longrightarrow  &  \left[ \p_r \, ^{s}U \right] &  \longrightarrow  & \mE\, ^{s}U^{\wedge}  && \\[1mm]
  &&\hspace{-8mm}{}_{- \uom}&\hspace{2mm}{} \hspace{11mm}{}_{\p_r}&\hspace{6mm}{}_{-\uom}&\hspace{11mm}{}_{-r}& \hspace{6mm}{}_{- \uom} &&\\[-2mm]
  && \uparrow & \diagdown \hspace{-1.3mm} \nearrow & \uparrow & \diagdown \hspace{-1.3mm} \nearrow & \uparrow &&\\[-1.1mm]
  && \downarrow & \diagup \hspace{-1.3mm} \searrow & \downarrow & \diagup \hspace{-1.3mm} \searrow& \downarrow  &&\\[-2mm]
  &&\hspace{-4.5mm}{}^{\uom}&\hspace{2mm}{} \hspace{12mm}{}^{-\p_r}&\hspace{5mm}{}^{\uom}&\hspace{2mm}{} \hspace{10mm}{}^{r}&\hspace{6mm}{}^{\uom}&&\\[1mm]
  && ^{s}U & \longrightarrow & \uom\, \left[\p_r\, ^{s}U \right] & \longrightarrow & \uom\, \mE\, ^{s}U^{\wedge}&&\\[-1mm]
  &&               & {}_{\uom \p_r}  &           & -\ux& &&
 \end{array}
$$
which can be compressed into
$$
\begin{array}{ccccccccc}
  &&    & {}_{\mE}   &                             &&&&\\
  &&  ^{s}U^{\wedge} = -\uom\,  ^{s}U & \longrightarrow  &  \mE\, ^{s}U^{\wedge} &&&& \\[1mm]
  &&\hspace{-8mm}{}_{- \uom}&&\hspace{6mm}{}_{-\uom}&&&&\\[-2mm]
  && \uparrow &   & \uparrow &&&&\\[-1.1mm]
  && \downarrow &   & \downarrow &&&&\\[-2mm]
  &&\hspace{-4.5mm}{}^{\uom}&&\hspace{5mm}{}^{\uom}&&&&\\[1mm]
  && ^{s}U & \longrightarrow & \uom\,  \mE\, ^{s}U^{\wedge} &&&&\\[-1mm]
  &&               & {}_{\mE}  &           &&&&
 \end{array}
$$
indeed confirming the signum--pair of operators $(\mE, \mE)$.
}
\end{remark}

\begin{remark}
{\rm
The commutative scheme of Definition \ref{signumradderiv} again shows the action of the operator $\dirac_{rad} = \uom\, \p_r$ on the signumdistribution $^{s}U$ resulting in the signumdistribution $\uom\, \p_r\, ^{s}U$ given by the equivalence class
$$
[\uom\, \p_r\, ^{s}U] = \uom\, \left[\p_r\, ^{s}U \right] = \uom\, \left[  -\, \frac{1}{\ux} \, \mE\, ^{s}U^{\wedge}  \right] = \uom\, \left[ \frac{1}{\ux} \, \mE\, \uom\, ^{s}U  \right] = \left[ -\, \frac{1}{\ux} \, \mE\, ^{s}U  \right] 
$$
confirming the signum--pair of operators $[\uom\, \p_r\, , \uom\, \p_r]$.\\
When, in particular, $^{s}U$ is a radial signumdistribution: $^{s}U = {^{s}}U^{rad}$, the action of the Dirac operator $\dirac$ on $^{s}U^{rad}$ reduces to
$$
\dirac\, ^{s}U^{rad} = \left[ \uom\, \p_r\, ^{s}U^{rad}  \right] = \left[  \frac{1}{\ux}\, \mE\, \uom\, ^{s}U^{rad} \right]
$$
which is in agreement with the action of $\dirac_{rad}$ on the signumdistribution $\uom\, T$, as contained in the commutative scheme of Definition \ref{radderiv}.
}
\end {remark}

\begin{definition}
\label{signumangderiv}
The angular $\p_{\uom}$--derivative of a scalar--valued signumdistribution $^{s}U$ is the distribution $\p_{\uom}\, ^{s}U = \p_{\uom}\, \uom\, ^{s}U^{\wedge}$ 
according to (the boldface part of) the commutative scheme
$$
\begin{array}{ccccccccc}
  &&    & {}_{\p_{\uom} \uom}   &                             &&&&\\
  && {\bf ^{s}U^{\wedge} = -\uom\,  ^{s}U } & \longrightarrow  & {\bf \p_{\uom} \, ^{s}U } &&&& \\[1mm]
  &&\hspace{-8mm}{}_{- \uom}&\hspace{2mm}{} \hspace{11mm}{}_{\p_{\uom}}&\hspace{6mm}{}_{-\uom}&&&&\\[-2mm]
  && \uparrow & \diagdown \hspace{-1.3mm} \nearrow & \uparrow &&&&\\[-1.1mm]
  && \downarrow & \diagup \hspace{-1.3mm} \searrow & \downarrow &&&&\\[-2mm]
  &&\hspace{-4.5mm}{}^{\uom}&\hspace{2mm}{} \hspace{12mm}{}^{\uom \p_{\uom} \uom}&\hspace{5mm}{}^{\uom}&&&&\\[1mm]
  &&  {\bf ^{s}U} & \longrightarrow & \uom\, \p_{\uom}  \, ^{s}U &&&&\\[-1mm]
  &&               & {}_{\uom \p_{\uom}}  &           &&&&
 \end{array}
$$
\end{definition}

\begin{remark}
{\rm
The commutative scheme of Definition \ref{signumangderiv} establishes the cross--pair of operators $(\uom\, \p_{\uom}\, \uom\, , \p_{\uom} )$ and confirms the signum--pair of operators $(-\, \p_{\uom}\, \uom\, , \, \Gamma) = ( -\, \Gamma + (m-1)\, {\bf 1}\, , \, \Gamma)$.
}
\end{remark}

  
\section{Composite actions of two operators}
\label{combined}
  

In the preceding sections we were able to define the actions on (signum)distributions of the operators $r$, $\uom$, $\p_r$, and $\p_{\uom}$. In Section \ref{intro} it was argued that the composite action by any two of those operators should lead to a {\em legal} action on distributions. Let us find out now if this is indeed the case.

\subsection{}

Multiplication of a distribution $T$ by the analytic function $r^2 = \sum_{j=1}^m\, x_j^2 = -\, \ux^2$ is well defined. Through the following commutative scheme it is shown that $r(rT) = r^2\, T$:
$$
\begin{array}{ccccccccc}
  &&    & {}_{-\, \ux}   &                             &  {}_{\ux} &&\\
  &&  T & \longrightarrow  &  -\, \ux\, T &  \longrightarrow  & r^2\, T  && \\[1mm]
  &&\hspace{-8mm}{}_{- \uom}&\hspace{2mm}{} \hspace{11mm}{}_{}&\hspace{6mm}{}_{-\uom}&\hspace{11mm}{}_{r}& \hspace{6mm}{}_{- \uom} &&\\[-2mm]
  && \uparrow & \diagdown \hspace{-1.3mm} \phantom{\nearrow} & \uparrow & \phantom{\diagdown} \hspace{-0.8mm} \nearrow & \uparrow &&\\[-1.1mm]
  && \downarrow & \phantom{\diagup} \hspace{-0.8mm} \searrow & \downarrow & \diagup \hspace{-1.3mm} \phantom{\searrow} & \downarrow  &&\\[-2mm]
  &&\hspace{-4.5mm}{}^{\uom}&\hspace{2mm}{} \hspace{12mm}{}^{r}&\hspace{5mm}{}^{\uom}&\hspace{2mm}{} \hspace{10mm}{}^{}&\hspace{6mm}{}^{\uom}&&\\[1mm]
  && \uom\, T & \longrightarrow & r\, T & \longrightarrow & \uom\, r^2 T &&\\[-1mm]
  &&               & {}_{-\, \ux}  &           &{}_{\ux}& &&
 \end{array}
$$

\subsection{}

Multiplication of a distribution $T$ by the analytic function $\ux = r\, \uom$ is well defined. Through the commutative scheme of Definition \ref{multiplicr} it is shown that $r(\uom\, T) = \ux\, T$.

\subsection{}

The action of the Euler operator $\mE = \sum_{j=1}^m\, x_j\, \p_{x_j}$ on a distribution is well defined. Through the following commutative scheme it is shown that $r(\p_r\, T) = \mE\, T$:
$$
\begin{array}{ccccccccc}
  &&    & {}_{-\, \uom\, \p_r}   &                             &  {}_{\ux} &&\\
  &&  T & \longrightarrow  &  [-\, \uom\, \p_r\, T] = \left[ \frac{1}{\ux}\, \mE\, T  \right] &  \longrightarrow  & \mE\, T  && \\[1mm]
  &&\hspace{-8mm}{}_{- \uom}&\hspace{2mm}{} \hspace{11mm}{}_{}&\hspace{6mm}{}_{-\uom}&\hspace{11mm}{}_{r}& \hspace{6mm}{}_{- \uom} &&\\[-2mm]
  && \uparrow & \diagdown \hspace{-1.3mm} \phantom{\nearrow} & \uparrow & \phantom{\diagdown} \hspace{-0.8mm} \nearrow & \uparrow &&\\[-1.1mm]
  && \downarrow & \phantom{\diagup} \hspace{-0.8mm} \searrow & \downarrow & \diagup \hspace{-1.3mm} \phantom{\searrow} & \downarrow  &&\\[-2mm]
  &&\hspace{-4.5mm}{}^{\uom}&\hspace{2mm}{} \hspace{12mm}{}^{\p_r}&\hspace{5mm}{}^{\uom}&\hspace{2mm}{} \hspace{10mm}{}^{}&\hspace{6mm}{}^{\uom}&&\\[1mm]
  && \uom\, T & \longrightarrow & [\p_r\, T] = \uom\, \left[ \frac{1}{\ux}\, \mE\, T  \right] & \longrightarrow & \uom\, \mE\, T &&\\[-1mm]
  &&               & {}_{-\, \uom\, \p_r}  &           &{}_{\ux}& &&
 \end{array}
$$
Moreover the signum--pair of operators $(\mE \, , \, \mE)$ is confirmed.

\subsection{}

The action of the operator $\ux\, \Gamma =  \ux\, \left( \sum_{j<k} \, e_j e_k (x_j \p{x_k} - x_k \p_{x_j}) \right)$ on a distribution is well defined. Through the following commutative scheme it is shown that $r(\p_{\uom}\, T) = \ux\, \Gamma\, T$:
$$
\begin{array}{ccccccccc}
  &&    & {}_{-\, \uom\, \p_{\uom}}   &                             &  {}_{\ux} &&\\
  &&  T & \longrightarrow  &  \Gamma\, T &  \longrightarrow  & \ux\, \Gamma\, T  && \\[1mm]
  &&\hspace{-8mm}{}_{- \uom}&\hspace{2mm}{} \hspace{11mm}{}_{}&\hspace{6mm}{}_{-\uom}&\hspace{11mm}{}_{r}& \hspace{6mm}{}_{- \uom} &&\\[-2mm]
  && \uparrow & \diagdown \hspace{-1.3mm} \phantom{\nearrow} & \uparrow & \phantom{\diagdown} \hspace{-0.8mm} \nearrow & \uparrow &&\\[-1.1mm]
  && \downarrow & \phantom{\diagup} \hspace{-0.8mm} \searrow & \downarrow & \diagup \hspace{-1.3mm} \phantom{\searrow} & \downarrow  &&\\[-2mm]
  &&\hspace{-4.5mm}{}^{\uom}&\hspace{2mm}{} \hspace{12mm}{}^{\p_{\uom}}&\hspace{5mm}{}^{\uom}&\hspace{2mm}{} \hspace{10mm}{}^{}&\hspace{6mm}{}^{\uom}&&\\[1mm]
  && \uom\, T & \longrightarrow & \p_{\uom}\, T & \longrightarrow & \ux\, \p_{\uom}\, T &&\\[-1mm]
  &&               & {}_{-\, \p_{\uom}\, \uom}  &           &{}_{\ux}& &&
 \end{array}
$$
and the signum--pair of operators $(r\, \p_{\uom} \, , \, -\, r\, \uom\, \p_{\uom}\, \uom) = (r\, \p_{\uom} \, , \, -\, r\, \p_{\uom} + (m-1)\, r\, \uom)$ is established.

\subsection{}

It is clear that $\uom\, (\uom\, T) = - T$.

\subsection{}

The action of the operator $\uom\, \p_r$ on a distribution is well defined, albeit not uniquely but through an equivalence class instead, see (\ref{equirad}). Definition \ref{radderiv} implies that
$\uom\, [\p_r\, T] = [(\uom\, \p_r)\, T]$.

\subsection{}

The action of the operator $\uom\, \p_{\uom} = -\, \Gamma$ on a distribution is well defined. Definition \ref{angderiv} implies that
$\uom\, (\p_{\uom}\, T) = -\, \Gamma\, T$.

\subsection{}

The action of the operator $\p_r^2$ on a distribution was defined in Section \ref{laplace} by the equivalence class
$$
\p_r^2\, T = \left[ -\, (\uom\, \p_r)^2\, T \right] = S_2 + \delta(\ux)\, c_2 - \sum_{j=1}^m\, c_{1,j}\, \p_{x_j} \delta(\ux)
$$
for arbitrary constants $c_2$ and $ c_{1,j}, j=1,\ldots,m$ and any distribution $S_2$ such that $\ux\, S_2 = \mE\, \underline{S}_1$ with $\ux\, \underline{S}_1 = -\, \mE\, T$. This is in complete agreement with the commutative scheme
$$
\begin{array}{ccccccccc}
  &&    & {}_{-\, \uom\, \p_r}   &                             &  {}_{\uom\, \p_r} &&\\
  &&  T & \longrightarrow  &  [-\, \uom\, \p_r\, T] = \left[ \frac{1}{\ux}\, \mE\, T  \right] &  \longrightarrow  & \left[ -\, (\uom\, \p_r)^2\, T \right]  && \\[1mm]
  &&\hspace{-8mm}{}_{- \uom}&\hspace{2mm}{} \hspace{11mm}{}_{}&\hspace{6mm}{}_{-\uom}&\hspace{11mm}{}_{\p_r}& \hspace{6mm}{}_{- \uom} &&\\[-2mm]
  && \uparrow & \diagdown \hspace{-1.3mm} \phantom{\nearrow} & \uparrow & \phantom{\diagdown} \hspace{-0.8mm} \nearrow & \uparrow &&\\[-1.1mm]
  && \downarrow & \phantom{\diagup} \hspace{-0.8mm} \searrow & \downarrow & \diagup \hspace{-1.3mm} \phantom{\searrow} & \downarrow  &&\\[-2mm]
  &&\hspace{-4.5mm}{}^{\uom}&\hspace{2mm}{} \hspace{12mm}{}^{\p_r}&\hspace{5mm}{}^{\uom}&\hspace{2mm}{} \hspace{10mm}{}^{}&\hspace{6mm}{}^{\uom}&&\\[1mm]
  && \uom\, T & \longrightarrow & [\p_r\, T] = \uom\, \left[ \frac{1}{\ux}\, \mE\, T  \right] & \longrightarrow & \uom\, \left[ -\, (\uom\, \p_r)^2\, T \right] &&\\[-1mm]
  &&               & {}_{-\, \uom\, \p_r}  &           &{}_{\uom\, \p_r}& &&
 \end{array}
$$
which moreover implies the signum--pair of operators $[\p_{r}^2 \, , \, \p_{r}^2]$.

\subsection{}
Start with the observation that for a distribution $T$, $\p_r\, \p_{\uom}\, T = \uom\, \p_r\, (-\, \uom\, \p_{\uom})\, T =  -\, \left[ \frac{1}{\ux}\, \mE\, \Gamma\, T  \right]$ to see that the action of the operator $\p_r\, \p_{\uom}$ on a distribution is well--defined, though not uniquely. Then the commutative scheme
$$
\begin{array}{ccccccccc}
  &&    & {}_{-\, \uom\, \p_{\uom}}   &                             &   {}_{ \uom\, \p_r}&&\\
  &&  T & \longrightarrow  &  \Gamma\, T &  \longrightarrow  & [\p_r\, \p_{\uom}\, T]  && \\[1mm]
  &&\hspace{-8mm}{}_{- \uom}&\hspace{2mm}{} \hspace{11mm}{}_{}&\hspace{6mm}{}_{-\uom}&\hspace{11mm}{}_{\p_{r}}& \hspace{6mm}{}_{- \uom} &&\\[-2mm]
  && \uparrow & \diagdown \hspace{-1.3mm} \phantom{\nearrow} & \uparrow & \phantom{\diagdown} \hspace{-0.8mm} \nearrow & \uparrow &&\\[-1.1mm]
  && \downarrow & \phantom{\diagup} \hspace{-0.8mm} \searrow & \downarrow & \diagup \hspace{-1.3mm} \phantom{\searrow} & \downarrow  &&\\[-2mm]
  &&\hspace{-4.5mm}{}^{\uom}&\hspace{2mm}{} \hspace{12mm}{}^{\p_{\uom}}&\hspace{5mm}{}^{\uom}&\hspace{2mm}{} \hspace{10mm}{}^{}&\hspace{6mm}{}^{\uom}&&\\[1mm]
  && \uom\, T & \longrightarrow & \p_{\uom}\, T & \longrightarrow & [-\, \p_r\, \Gamma\, T]  &&\\[-1mm]
  &&               & {}_{-\, \p_{\uom}\, \uom}  &           &{}_{ \uom\, \p_r} & &&
 \end{array}
$$
shows that indeed $\p_r\, (\p_{\uom}\, T) = \p_r\, \p_{\uom}\, T$.

\subsection{}
We know that the action of the operator $\p_{\uom}^2$ on a distribution is uniquely defined. Applying twice the commutative scheme of Definition \ref{angderiv} we obtain
$$
\begin{array}{ccccccccc}
  &&    & {}_{-\, \uom\, \p_{\uom}}   &                             &   {}_{ \p_{\uom}\, \uom}&&\\
  &&  T & \longrightarrow  &  \Gamma\, T &  \longrightarrow  & \p_{\uom}^2\, T  && \\[1mm]
  &&\hspace{-8mm}{}_{- \uom}&\hspace{2mm}{} \hspace{11mm}{}_{}&\hspace{6mm}{}_{-\uom}&\hspace{11mm}{}_{\p_{\uom}}& \hspace{6mm}{}_{- \uom} &&\\[-2mm]
  && \uparrow & \diagdown \hspace{-1.3mm} \phantom{\nearrow} & \uparrow & \phantom{\diagdown} \hspace{-0.8mm} \nearrow & \uparrow &&\\[-1.1mm]
  && \downarrow & \phantom{\diagup} \hspace{-0.8mm} \searrow & \downarrow & \diagup \hspace{-1.3mm} \phantom{\searrow} & \downarrow  &&\\[-2mm]
  &&\hspace{-4.5mm}{}^{\uom}&\hspace{2mm}{} \hspace{12mm}{}^{\p_{\uom}}&\hspace{5mm}{}^{\uom}&\hspace{2mm}{} \hspace{10mm}{}^{}&\hspace{6mm}{}^{\uom}&&\\[1mm]
  && \uom\, T & \longrightarrow & \p_{\uom}\, T & \longrightarrow & \uom\,  \p_{\uom}^2\, T  &&\\[-1mm]
  &&               & {}_{-\, \p_{\uom}\, \uom}  &           &{}_{ \uom\, \p_{\uom}} & &&
 \end{array}
$$
showing that indeed $\p_{\uom}\, ( \p_{\uom}\, T) = \p_{\uom}^2\, T$ and confirming the signum--pair of operators $(\p_{\uom}^2 \, , \, \p_{\uom}^2)$.

  
\section{Division of (signum)distributions by $r$}
\label{division}
  

Division of a standard distribution $T$ by an analytic function $\alpha(\ux)$ resulting in an equivalent class of distributions $S$ such that $\alpha(\ux)\, S = T$, we expect the division of a standard distribution by the non--analytic function $r$ to lead to an equivalence class of signumdistributions. Let us make this precise.

\begin{definition}
\label{T/r}
The quotient of a scalar distribution $T$ by the radial distance $r$ is the equivalece class of signumdistributions
$$
\left[ \frac{1}{r}\, T \right] = \uom\, \left[ \frac{1}{\ux}\, T  \right] = \uom\, (\underline{S} + \delta(\ux)\, \underline{c} ) = \uom\, \underline{S} + \uom\, \delta(\ux)\, \underline{c} = \underline{S}^{\vee} + \delta(\ux)^{\vee}\, \underline{c}
$$
for any vector-valued distribution $\underline{S}$ for which $\ux\, \underline{S} = T$,
according to (the boldface part of) the commutative scheme
$$
\begin{array}{ccccccccc}
  &&    & {}_{\frac{1}{\ux}}   &                             &  &&&\\
  && {\bf T} & \longrightarrow  & {\bf \left[ \frac{1}{\ux}\, T \right]} & &&& \\[1mm]
  &&\hspace{-8mm}{}_{- \uom} &\hspace{2mm}{} \hspace{11mm}{}_{- \frac{1}{r}}&\hspace{6mm}{}_{-\uom}&&&&\\[-2mm]
  && \uparrow & \diagdown \hspace{-1.3mm} \nearrow & \uparrow &&&&\\[-1.1mm]
  && \downarrow & \diagup \hspace{-1.3mm} \searrow & \downarrow &&&&\\[-2mm]
  &&\hspace{-4.5mm}{}^{\uom}&\hspace{2mm}{} \hspace{10mm}{}^{\frac{1}{r}}&\hspace{5mm}{}^{\uom}&&&&\\[1mm]
  && T^{\vee} = \uom\, T & \longrightarrow & {\bf \left[ \frac{1}{r}\, T \right]} &&&&\\[-1mm]
  &&               & {}_{\frac{1}{\ux}}  &           &&&&
 \end{array}
$$
\end{definition}

\begin{example}
{\em
Let us illustrate Definition \ref{T/r} by the case of the delta-distribution: $T = \delta(\ux)$. As $\ux\, \dirac\, \delta(\ux) = m\, \delta(\ux)$ and $\ux\, \delta(\ux) = 0$ we have
$$
\frac{1}{\ux}\, \delta(\ux) = \frac{1}{m}\, \dirac\, \delta(\ux) + \delta(\ux)\, \underline{c}_0
$$
with $\underline{c}_0$ an arbitrary constant vector. It then follows that
$$
\left[ \frac{1}{r}\, \delta(\ux) \right] = \uom\, \left[ \frac{1}{\ux}\, \delta(\ux)  \right] = \uom\, \left[ \frac{1}{m}\, \dirac\, \delta(\ux)  + \delta(\ux)\, \underline{c}_0 \right] = \frac{1}{m}\, \uom\, \dirac\, \delta(\ux) + \uom\, \delta(\ux)\, \underline{c}_0 = \frac{1}{m}\, (\dirac\, \delta(\ux))^{\vee} + \delta(\ux)^{\vee}\, \underline{c}_0
$$
or, in view of the definition of $\p_r\, \delta(\ux)$,
$$
\left[ \invr\, \delta(\ux) \right] = -\, \frac{1}{m}\, \p_r\, \delta(\ux) + \uom\, \delta(\ux)\, \underline{c}_0
$$
However in this particular case of the delta-distribution  it turns out that $\invr\, \delta(\ux)$ is uniquely determined. Indeed, as $\p_r\, \delta(\ux)$ is a radial signumdistribution and as we expect the signumdistribution $\invr\, \delta(\ux)$ to be SO$(m)$--invariant as well, the arbitrary vector constant $\underline{c}_0$ should be zero, eventually leading to the formulae
$$
\invr\, \delta(\ux)  = -\, \frac{1}{m}\, \p_r\, \delta(\ux) \quad , \quad \invr\, \uom\,  \delta(\ux)  =  -\, \frac{1}{m}\, \dirac\, \delta(\ux)
$$
and
$$
\invux\, \delta(\ux)  = \frac{1}{m}\, \dirac\, \delta(\ux) \quad , \quad \invux\, \uom\, \delta(\ux)  = -\, \frac{1}{m}\, \p_r\, \delta(\ux)
$$
which fit into the commutative scheme
$$
\begin{array}{ccccccccc}
  &&    & {}_{\frac{1}{\ux}}   &                             &  &&&\\
  && \delta(\ux) & \longrightarrow  & \frac{1}{m}\, \dirac \delta(\ux) & &&& \\[1mm]
  &&\hspace{-8mm}{}_{- \uom} &\hspace{2mm}{} \hspace{11mm}{}_{- \frac{1}{r}}&\hspace{6mm}{}_{-\uom}&&&&\\[-2mm]
  && \uparrow & \diagdown \hspace{-1.3mm} \nearrow & \uparrow &&&&\\[-1.1mm]
  && \downarrow & \diagup \hspace{-1.3mm} \searrow & \downarrow &&&&\\[-2mm]
  &&\hspace{-4.5mm}{}^{\uom}&\hspace{2mm}{} \hspace{10mm}{}^{\frac{1}{r}}&\hspace{5mm}{}^{\uom}&&&&\\[1mm]
  && \delta(\ux)^{\vee} = \uom\, \delta(\ux) & \longrightarrow & -\, \frac{1}{m}\, \p_r\, \delta(\ux) &&&&\\[-1mm]
  &&               & {}_{\frac{1}{\ux}}  &           &&&&
 \end{array}
$$
For the more general case of the division of the delta-distribution by natural powers of $r$ we refer to \cite{bsv}.
}
\end{example}

\begin{remark}
{\rm
The observation made in the above example about the unique definition of $\invr\, \delta(\ux)$ holds more generally for all radial distributions, as it might be expected that, for  a given radial distribution $T^{rad}$, $\invr\, T^{rad}$ is a radial signumdistribution as well. If $[\invux\, T^{rad}] = \underline{S}$ with $\ux\, \underline{S} = T^{rad}$, then $\invr\, T^{rad} = \uom\, \underline{S}$ should be a radial signumdistribution. We claim that such a distribution $\underline{S}$ is unique. Indeed, if $\ux \underline{S}_1 = \ux \underline{S}_2 = T^{rad}$ then $\ux\, (\underline{S}_1 - \underline{S}_2) = 0$, whence $\underline{S}_1 - \underline{S}_2 = \delta(\ux)\, \underline{c}$. As $\uom\, (\underline{S}_1 - \underline{S}_2) = \uom\, \delta(\ux)\, \underline{c} $ should be a radial signumdistribution, the arbitrary constant $\underline{c}$ should be zero. The conclusion is that $\invr\, T^{rad}$ is uniquely defined by the {\em radial} signumdistribution $\uom\, \underline{S}$ for which $\ux\, \underline{S} = T^{rad}$. It then also holds that 
$
\invux\, T^{rad} = \underline{S}
$ is uniquely defined.
}
\end{remark}

\begin{remark}
{\rm
The commutative scheme of Definition \ref{T/r} involves the signum--pair of operators $[\frac{1}{\ux}\, , \frac{1}{\ux}]$ and the cross--pair of operators $[\invr\, , -\, \invr]$. The signum--pair of operators $[\invrsq\, , \, \invrsq]$ follows at once.
}
\end{remark}

It is also interesting and useful to make the division by $r$ of a signumdistribution explicit, because it will lead to the definition of the action of the angular part $\dirac_{ang} = \frac{1}{r}\, \p_{\uom}$ of the Dirac operator on a signumdistribution, leading in its turn to the definition of the action of the Dirac operator on a signumdistribution.

\begin{definition}
\label{U/r}
The quotient of a scalar-valued signumdistribution $^{s}U$ by the radial distance $r$ is the equivalece class of distributions
$$
\left[ \frac{1}{r}\, ^{s}U \right] =  \left[ \frac{1}{\ux}\, \uom\, ^{s}U  \right] = S + \delta(\ux)\, c 
$$
for any scalar-valued distribution $S$ for which $\ux\, S = \uom\, ^{s}U$,
according to (the boldface part of) the commutative scheme
$$
\begin{array}{ccccccccc}
  &&    & {}_{-\frac{1}{\ux}}   &                             &  &&&\\
  && {\bf -\uom ^{s}U} & \longrightarrow  & {\bf \left[ \frac{1}{r}\, ^{s}U \right]} & &&& \\[1mm]
  &&\hspace{-8mm}{}_{- \uom} &\hspace{2mm}{} \hspace{11mm}{}_{ \frac{1}{r}}&\hspace{6mm}{}_{-\uom}&&&&\\[-2mm]
  && \uparrow & \diagdown \hspace{-1.3mm} \nearrow & \uparrow &&&&\\[-1.1mm]
  && \downarrow & \diagup \hspace{-1.3mm} \searrow & \downarrow &&&&\\[-2mm]
  &&\hspace{-4.5mm}{}^{\uom}&\hspace{2mm}{} \hspace{10mm}{}^{- \frac{1}{r}}&\hspace{5mm}{}^{\uom}&&&&\\[1mm]
  && {\bf ^{s}U} & \longrightarrow &  \uom\, \left[ \frac{1}{\ux}\, \uom\, ^{s}U \right] &&&&\\[-1mm]
  &&               & {}_{-\frac{1}{\ux}}  &           &&&&
 \end{array}
$$
\end{definition}

\begin{remark}
{\rm
The commutative scheme of Definition \ref{U/r} confirms the signum--pair of operators $[-\, \frac{1}{\ux}\, , -\,\frac{1}{\ux}]$ and the cross--pair of operators $[-\, \invr\, ,  \invr]$. 
}
\end{remark}

As we know how to act with the operator $\p_{\uom}$ on a distribution (see Definition \ref{angderiv}) and how to act with the operator $\invr$ on a signumdistribution (see Definition \ref{U/r}) we are now in the position to check the action on a distribution of the composition of both operators, viz. the angular part $\dirac_{ang}$ of the Dirac operator. The outcome should match result (\ref{equiang}), which indeed it does seen the following commutative scheme 
$$
\begin{array}{ccccccccc}
  &&    & {}_{-\, \uom \p_{\uom}}   &                             &  {}_{-\, \frac{1}{\ux}} &&\\
  &&  T & \longrightarrow  &  \Gamma\, T &  \longrightarrow  & \left[ \invr\, \p_{\uom}\, T \right]  = \left[ -\, \frac{1}{\ux}\, \Gamma\, T\right]  && \\[1mm]
  &&\hspace{-8mm}{}_{- \uom}&\hspace{2mm}{} \hspace{11mm}{}_{}&\hspace{6mm}{}_{-\uom}&\hspace{11mm}{}_{\invr}& \hspace{6mm}{}_{- \uom} &&\\[-2mm]
  && \uparrow & \diagdown \hspace{-1.3mm} \phantom{\nearrow} & \uparrow & \phantom{\diagdown} \hspace{-0.8mm} \nearrow & \uparrow &&\\[-1.1mm]
  && \downarrow & \phantom{\diagup} \hspace{-0.8mm} \searrow & \downarrow & \diagup \hspace{-1.3mm} \phantom{\searrow} & \downarrow  &&\\[-2mm]
  &&\hspace{-4.5mm}{}^{\uom}&\hspace{2mm}{} \hspace{12mm}{}^{\p_{\uom}}&\hspace{5mm}{}^{\uom}&\hspace{2mm}{} \hspace{10mm}{}^{}&\hspace{6mm}{}^{\uom}&&\\[1mm]
  && \uom\, T & \longrightarrow & \p_{\uom}\, T & \longrightarrow & \uom\, \left[ -\, \frac{1}{\ux}\, \Gamma\, T\right] &&\\[-1mm]
  &&               & {}_{-\, \p_{\uom}\, \uom}  &           &{}_{-\, \frac{1}{\ux}}& &&
 \end{array}
$$
which moreover confirms the signum--pair of operators $$[\invr\, \p_{\uom} \, , \, -\, \invr\, \uom\, \p_{\uom}\, \uom ] = [\invr\, \p_{\uom} \, , \, -\, \invr\, \p_{\uom} + (m-1)\, \invr\, \uom ]$$
In the same order of ideas we can define the action of $\dirac_{ang} = \invr\, \p_{\uom}$ on a signumdistribution through the commutative scheme
$$
\begin{array}{ccccccccc}
  &&    & {}_{\p_{\uom}\, \uom}   &                             &  \frac{1}{\ux} &&\\
  &&  ^{s}U^{\wedge} = -\uom\,  ^{s}U & \longrightarrow  &   \p_{\uom} \, ^{s}U  &  \longrightarrow  & \left[ \frac{1}{\ux}\, \p_{\uom}\, ^{s}U \right] && \\[1mm]
  &&\hspace{-8mm}{}_{- \uom}&\hspace{2mm}{} \hspace{11mm}{}_{\p_{\uom}}&\hspace{6mm}{}_{-\uom}&\hspace{11mm}{}_{}& \hspace{6mm}{}_{- \uom} &&\\[-2mm]
  && \uparrow & \phantom{\diagdown} \hspace{-0.5mm} \nearrow & \uparrow & \diagdown \hspace{-1.3mm} \phantom{\nearrow} & \uparrow &&\\[-1.1mm]
  && \downarrow & \diagup \hspace{-1.3mm} \phantom{\searrow} & \downarrow & \phantom{\diagup} \hspace{-0.6mm} \searrow& \downarrow  &&\\[-2mm]
  &&\hspace{-4.5mm}{}^{\uom}&\hspace{2mm}{} \hspace{12mm}{}^{}&\hspace{5mm}{}^{\uom}&\hspace{2mm}{} \hspace{8mm}{}^{\invr}&\hspace{6mm}{}^{\uom}&&\\[1mm]
  && ^{s}U & \longrightarrow & \uom\, \p_{\uom}\, ^{s}U & \longrightarrow & \left[ \frac{1}{r}\, \p_{\uom}\, ^{s}U \right]&&\\[-1mm]
  &&               & {}_{\uom \p_{\uom}}  &           & \frac{1}{\ux} & &&
 \end{array}
$$
in other words
$$
\left[ \dirac_{ang}\, ^{s}U  \right] = \left[ \invr\, \p_{\uom}\, ^{s}U  \right] = \uom\, \left[ \frac{1}{\ux}\,   \p_{\uom}\, ^{s}U \right] = \left[ -\,  \frac{1}{\ux}\, \Gamma\,   ^{s}U \right]
$$
which leads to the signum--pair of operators $[\invr\, (-\uom)\, \p_{\uom}\, \uom\, , \invr\, \p_{\uom}]$ and the cross--pair of operators $[\invr\, \p_{\uom}\, \uom\, , -\, \invr\, \uom\, \p_{\uom}]$.\\

Combining the actions on a signumdistribution of the radial and angular parts of the Dirac operator, we are able to define the action of the Dirac operator itself on a signumdistribution.

\begin{definition}
\label{diracsignumdistrib}
The action of the Dirac operator $\dirac$ on the signumdsitribution $^{s}U$ is given by the equivalence class of signumdistributions
\begin{eqnarray*}
\left[  \dirac\, ^{s}U \right] &=&  \left[ (\uom\, \p_r + \invr\, \p_{\uom} )\, ^{s}U \right]\\
&=& \left[ -\, \frac{1}{\ux}\, \mE\,   ^{s}U \right] + \left[ -\, \frac{1}{\ux}\, \Gamma\,   ^{s}U \right]\\
&=& \left[  -\, \frac{1}{\ux}\, (\mE + \Gamma)\,  ^{s}U \right]\\
&=& \left[  \frac{1}{\ux}\, (\ux\, \dirac)\, ^{s}U \right]
\end{eqnarray*}
according to the commutative scheme
$$
\begin{array}{ccccccccc}
  &&    & {}_{\underline{D}}   &                             &&&&\\
  &&  ^{s}U^{\wedge} = -\uom\,  ^{s}U & \longrightarrow  & \left[  \p_r\, ^{s}U + \frac{1}{\ux}\, \p_{\uom}\, ^{s}U \right] &&&& \\[1mm]
  &&\hspace{-8mm}{}_{- \uom}&\hspace{2mm}{} \hspace{11mm}{}_{}&\hspace{6mm}{}_{-\uom}&&&&\\[-2mm]
  && \uparrow & \phantom{\diagdown} \hspace{-1.3mm} \phantom{\nearrow} & \uparrow &&&&\\[-1.1mm]
  && \downarrow & \phantom{\diagup} \hspace{-1.3mm} \phantom{\searrow} & \downarrow &&&&\\[-2mm]
  &&\hspace{-4.5mm}{}^{\uom}&\hspace{2mm}{} \hspace{12mm}{}^{}&\hspace{5mm}{}^{\uom}&&&&\\[1mm]
  &&   ^{s}U & \longrightarrow & \left[ \dirac\, ^{s}U \right] &&&&\\[-1mm]
  &&               & {}_{\dirac}  &           &&&&
 \end{array}
$$
where we recall the operator $\underline{D}$ to be 
\begin{eqnarray*}
\underline{D} &=& \uom\, \p_r + \frac{1}{\ux}\, \p_{\uom}\, \uom\\
&=& \uom\, \p_r - \invr\, \uom\, \p_{\uom}\, \uom\\
&=& \uom\, \p_r - \invr\, \p_{\uom}\, + (m-1)\invr\, \uom
\end{eqnarray*}
\end{definition}

\begin{remark}
{\rm
The commutative scheme of Definition \ref{diracsignumdistrib} leads to the signum--pair of operators $[\underline{D}\, , \dirac]$, which in its turn leads to the signum--pair of operators $[{\bf Z}\, , \, \Delta]$. In this way the action of the Laplace operator $\Delta$ on signumdistributions is defined through the action of the operator ${\bf Z}$ on distributions; the result of this action is an equivalence class of signumdistributions.
Notice that while the signum--pairs of operators $(\dirac\, ,\, \underline{D})$ and $(\Delta\, , \, {\bf Z})$ are uniquely defined, the action result of the signum--pairs of operators  $[\underline{D}\, ,\, \dirac]$ and $[{\bf Z}\, , \, \Delta]$ are equivalent classes of (signum)distributions.

}
\end{remark}

\begin{example}
{\rm
Let us illustrate Definition \ref{diracsignumdistrib} with the following simple example. Consider the regular signumdistribition $\ux$ defined by
$$
\langle \ \ux \ , \ \uom\, \varphi \ \rangle = \langle \ \ux \, \uom \ , \  \varphi \ \rangle = \langle \ -\, r \ , \  \varphi \ \rangle = \int_{\mR^m}\, r\, \widetilde{\varphi}(r,\uom)\, d\ux
$$
for which
$$
\mE\, \ux = \ux
$$
and
$$
\Gamma\, \ux = (m-1)\ux
$$
whence
$$
\left[ \dirac\, \ux \right] = \left[ -\, \frac{1}{\ux}\, m\, \ux   \right] = \left[ -\, m \right] = -\, m + \delta(\ux)\, c
$$
As $\ux^{\wedge} = r$ and $\underline{D}\, r = \left[ m\, \uom \right]$, this result fits into the following commutative scheme:
$$
\begin{array}{ccccccccc}
  &&    & {}_{\underline{D}}   &                             &&&&\\
  &&  r & \longrightarrow  & \left[  m\, \uom \right] &&&& \\[1mm]
  &&\hspace{-8mm}{}_{- \uom}&\hspace{2mm}{} \hspace{11mm}{}_{}&\hspace{6mm}{}_{-\uom}&&&&\\[-2mm]
  && \uparrow & \phantom{\diagdown} \hspace{-1.3mm} \phantom{\nearrow} & \uparrow &&&&\\[-1.1mm]
  && \downarrow & \phantom{\diagup} \hspace{-1.3mm} \phantom{\searrow} & \downarrow &&&&\\[-2mm]
  &&\hspace{-4.5mm}{}^{\uom}&\hspace{2mm}{} \hspace{12mm}{}^{}&\hspace{5mm}{}^{\uom}&&&&\\[1mm]
  &&   \ux & \longrightarrow & \left[ -\, m \right] &&&&\\[-1mm]
  &&               & {}_{\dirac}  &           &&&&
 \end{array}
$$
}
\end{example}

Now as we know how to act with the multiplication operator $\invr$ on a signumdistribution, we can check the action on a distribution $T$ of the composite operator $(\invr\, \circ \, \p_r)\, T = \invr\, (\p_r\, T)$ which should coincide with the action $(\invr\, \p_r)\, T$, defined, though not uniquely, in Proposition \ref{laplaceparts} by
$\invr\, \p_r\, T = S_3 + \frac{1}{m}\, \sum_{j=1}^m\, c_{1,j}\, \p_{x_j} \delta(\ux) + c_3\, \delta(\ux)$ for arbitrary constant $c_3$ and any distribution $S_3$ such that $\ux\, S_3 = \underline{S}_1$ with $\ux\, \underline{S}_1 = -\, \mE\, T$. That this is indeed the case is shown by the following commutative scheme:
$$
\begin{array}{ccccccccc}
  &&    & {}_{-\, \uom\, \p_r}   &                             &  {}_{-\, \frac{1}{\ux}} &&\\
  &&  T & \longrightarrow  &  [-\, \uom\, \p_r\, T] = \left[ \frac{1}{\ux}\, \mE\, T  \right] &  \longrightarrow  & \left[ -\, \frac{1}{\ux^2}\, \mE\,T \right] =  \left[ \frac{1}{r^2}\, \mE\,T \right]  && \\[1mm]
  &&\hspace{-8mm}{}_{- \uom}&\hspace{2mm}{} \hspace{11mm}{}_{}&\hspace{6mm}{}_{-\uom}&\hspace{11mm}{}_{\invr}& \hspace{6mm}{}_{- \uom} &&\\[-2mm]
  && \uparrow & \diagdown \hspace{-1.3mm} \phantom{\nearrow} & \uparrow & \phantom{\diagdown} \hspace{-0.8mm} \nearrow & \uparrow &&\\[-1.1mm]
  && \downarrow & \phantom{\diagup} \hspace{-0.8mm} \searrow & \downarrow & \diagup \hspace{-1.3mm} \phantom{\searrow} & \downarrow  &&\\[-2mm]
  &&\hspace{-4.5mm}{}^{\uom}&\hspace{2mm}{} \hspace{12mm}{}^{\p_r}&\hspace{5mm}{}^{\uom}&\hspace{2mm}{} \hspace{10mm}{}^{}&\hspace{6mm}{}^{\uom}&&\\[1mm]
  && \uom\, T & \longrightarrow & [\p_r\, T] = \uom\, \left[ \frac{1}{\ux}\, \mE\, T  \right] & \longrightarrow & \uom\, \left[ \frac{1}{r^2}\, \mE\, T \right] &&\\[-1mm]
  &&               & {}_{-\, \uom\, \p_r}  &           &{}_{-\, \frac{1}{\ux}} & &&
 \end{array}
$$
implying the signum--pair of operators $[\invr\, \p_r \, , \, \invr\, \p_r]$.\\

Finally notice that, invoking the signum--pairs of operators $[ \p_r^2 \, , \,  \p_r^2]$, $[\invr\, \p_r \, , \, \invr\, \p_r]$ and $[\invrsq\, {\bf Z}^* \, , \, \invrsq\, \Delta^*]$, the signum--pair of operators $[{\bf Z}\, , \, \Delta]$ is confirmed.


\section{Cartesian derivatives of signumdistributions}
\label{signumcartderiv}


Recalling the signum--pairs of operators $(\dirac\, , \underline{D})$ and $[\underline{D}\, , \dirac]$, we expect the cartesian derivatives of (signum)distributions to show up in signum--pairs of operators of the form $(\p_{x_j}\, , d_j)$ and $[d_j\, , \p_{x_j}]$, $j=1,\ldots,m$, where the differential operators $d_j$ are related to the operator $\underline{D}$ in a similar way as the partial derivative operators $\p_{x_j}$ are related to the Dirac operator $\dirac$. Until now the operator $\uD$ was known only in terms of spherical co-ordinates. Let us show that $\uD$ indeed can be expressed in terms of cartesian co-ordinates. Taking into a	account that
\begin{eqnarray*}
\uom\, \p_r &=& -\, \frac{1}{\ux}\, \mE = \frac{\ux}{r^2}\, \mE\\
\invr\, \p_{\uom} &=&  -\, \frac{1}{\ux}\, \Gamma = \frac{\ux}{r^2}\, \Gamma\\
\invr\, \uom &=&  \frac{\ux}{r^2}
\end{eqnarray*}
we obtain, next to the well--known expression for the Dirac operator
$$
\dirac = \uom\, \p_r + \invr\, \pÐ{\uom} = -\, \invux\, (\mE + \Gamma) = \invrsq\, \ux\, (\mE + \Gamma)
$$
also the expression aimed at for the operator $\uD$:
\begin{eqnarray*}
\uD &=&  \uom\, \p_r - \invr\, \p_{\uom} + (m-1)\, \invr\, \uom\\
&=& \invrsq\, \ux\, (\mE - \Gamma + (m-1))\\
&=& \invrsq\, \ux\, (2\, \mE + (m-1)) - \dirac
\end{eqnarray*}

For $j=1,\ldots,m$ we define the operator $d_j$ as being the signum--partner of the cartesian derivative $\p_{x_j}$, in other words $(\p_{x_j}\, , d_j)$ is a signum--pair of operators. It thus holds for $j=1,\ldots,m$ that
\begin{eqnarray*}
d_j &=& \uom\, \p_{x_j}\, (-\uom)\\
&=& -\, \uom\, (\invr\, e_j - \ux\ \frac{x_j}{r^3}) + \p_{x_j}\\
&=& -\, \invr\, \uom\, e_j -\, \invr\, \omega_j + \p_{x_j}  = \invr\, e_j\, \uom + \invr\, \omega_j + \p_{x_j}\\
&=& -\, \invrsq\, \ux\, e_j - \invrsq\, x_j + \p_{x_j} = \invrsq\, e_j\, \ux + \invrsq\, x_j + \p_{x_j}
\end{eqnarray*}
Notice that the operator $d_j, j=1,\ldots,m$ contains the scalar part $\p_{x_j}$ and the bivector part $\invrsq\, (e_j\, \ux +x_j) = \invrsq\, e_j\, (e_1\, x_1 + \cdots + (e_j\, x_j)^{\wedge} +  \cdots + e_m\, x_m)$, as it should since it is the product of two vector operators. This bivector part is special in this sense that its multiplication by $\uom$ results into a vector operator as it should too.\\
As the operator $d_j$ clearly is a well--defined, but not uniquely defined, operator acting on distributions, also the signum--pairs of operators $[d_j\, , \p_{x_j}]$, $j=1,\ldots,m$ hold, enabling the definition of the cartesian derivative $\p_{x_j}$ of a signumdistribution as an equivalence class of signumdistributions, as is seen by the following commutative scheme:
$$
\begin{array}{ccccccccc}
  &&    & {}_{d_j}   &                             &&&&\\
  &&  ^{s}U^{\wedge} = -\uom\,  ^{s}U & \longrightarrow  & [  d_j\, (-\uom\, ^{s}U) ] &&&& \\[1mm]
  &&\hspace{-8mm}{}_{- \uom}&\hspace{2mm}{} \hspace{11mm}{}_{}&\hspace{6mm}{}_{-\uom}&&&&\\[-2mm]
  && \uparrow & \phantom{\diagdown} \hspace{-1.3mm} \phantom{\nearrow} & \uparrow &&&&\\[-1.1mm]
  && \downarrow & \phantom{\diagup} \hspace{-1.3mm} \phantom{\searrow} & \downarrow &&&&\\[-2mm]
  &&\hspace{-4.5mm}{}^{\uom}&\hspace{2mm}{} \hspace{12mm}{}^{}&\hspace{5mm}{}^{\uom}&&&&\\[1mm]
  &&   ^{s}U & \longrightarrow & \left[ \p_{x_j}\, ^{s}U \right] &&&&\\[-1mm]
  &&               & {}_{\p_{x_j}}  &           &&&&
 \end{array}
$$
\begin{example}
{\rm
Consider the regular distribution $\uom$ for which $\p_{x_j}\, \uom =  \invr\, e_j - \invr\, \omega_j\, \uom$, $j=1,\ldots,m$. It follows that $d_j\, ^{s}1 =  \uom\, (-\, \invr\, e_j + \invr\, \omega_j\, \uom) = -\, \invr\, \uom\, e_j - \invr\, \omega_j$. When considering the same locally integrable function $\uom$ as a regular signumdistribution we find $\p_{x_j}\, \uom = [\uom\, d_j\, (-\uom)\, \uom] = [\uom\, d_j\, 1] = [\uom\, (-\, \invr\, \uom\, e_j -\, \invr\, \omega_j + \p_{x_j} )\, 1] = [\invr\, e_j -\, \invr\, \omega_j\, \uom], j=1,\ldots,m$, this time an equivalent class of signumdistributions. }
\end{example}

\begin{remark}
{\rm
One may wonder why a cartesian derivative of a signumdistribution turns out to be an equivalent class of signumdistributions, and not merely a signumdistribution, keeping in mind that a cartesian derivative of a distribution is simply a distribution given by the following differentiation rule:
$$
\langle \, \p_{x_j}\, T \, , \, \varphi \, \rangle = -\, \langle \,  T \, , \, \p_{x_j}\,\varphi \, \rangle
$$
which after all is nothing else but a subtle form of ''partial integration`` based upon the product rule for differentiation
$$
(\p_{x_j}T)\, \varphi = \p_{x_j}\, (T\, \varphi) - T\, (\p_{x_j}\varphi)
$$
where $T$ is a regular distribution associated to a $C_1$--function.\\
And one may also wonder if it is possible to define cartesian differentiation of a signumdistribution directly through a similar formula, without having to invoke an action on the associated distribution as is done through the signum--pair of operators $[\, d_j\, , \, \p_{x_j} \,]$.\\
Let us find out by considering a scalar--valued $C_1$--function $f$ acting as a regular signumdistribution. For each scalar--valued test function $\varphi$ it holds that
$$
\p_{x_j}\, (f\, \uom\, \varphi) = (\p_{x_j}f)\, \uom\, \varphi + f\, (\p_{x_j}\uom)\, \varphi + f\, \uom\, (\p_{x_j}\varphi)
$$
and so
$$
\p_{x_j}\, (f\, \uom\, \varphi) = (\p_{x_j}f)\, \uom\, \varphi + f\, \uom\, (-\uom\, \invr\, e_j - \invr\, \omega_j)\, \varphi + f\, \uom\, (\p_{x_j}\varphi)
$$
whence
$$
\int_{\mR^m}\, \p_{x_j}\, (f\, \uom\, \varphi)\, d\ux = \int_{\mR^m}\, (\p_{x_j}f)\, \uom\, \varphi\, d\ux + \int_{\mR^m}\, f\, \uom\, (-\uom\, \invr\, e_j - \invr\, \omega_j + \p_{x_j})\,\varphi\, d\ux
$$
and finally
$$
\int_{\mR^m}\, (\p_{x_j}f)\, \uom\, \varphi\, d\ux \  = \ -\, \int_{\mR^m}\, f\, \uom\, d_j\varphi\, d\ux
$$
or, in the language of (signum)distributions,
$$
\langle\, \p_{x_j}f\, ,\,  \uom\, \varphi\, \rangle  = -\, \langle\, f\, ,\,  \uom\, (d_j\varphi)\, \rangle
$$
This result shows that when taking the $\p_{x_j}$--derivative of the signumdistribution $f$, the derivative shifts to the test function but then in the form of  the operator $d_j$. But there is more: the right--hand side of the last expression may also be written as
$$ -\, \langle\, f\, ,\,  \uom\, (\invrsq\, e_j\, \ux + \invrsq\, x_j + \p_{x_j})\, \varphi\, \rangle =  -\, \langle\, \invrsq\, f\, ,\,  \uom\, ( e_j\, \ux +  x_j )\, \varphi\, \rangle\\ -\, \langle\, f\, ,\,  \uom\, (\p_{x_j}\, \varphi)\, \rangle$$ invoking the test functions $( e_j\, \ux +  x_j )\, \varphi$ and $\p_{x_j}\, \varphi$\,.
Obviously the first term at the right--hand side of this last expression is indeed an equivalence class of signumdistributions seen the division of $f$ by the analytic function $r^2$ showing a zero at the origin.
}
\end{remark}

\begin{remark}
{\rm
As expected there is a relationship between the operators $d_j, j=1,\ldots,m$, $\dirac$ and $\uD$. A straightforward calculation shows that
\begin{eqnarray*}
\dirac &=& \phantom{-}\, \sum_{j=1}^m\, e_j\, d_j + (m-1)\, \invrsq\, \ux\\
\uD &=& -\, \sum_{j=1}^m\, e_j\, d_j  + 2\, \invrsq\, \ux\, \mE
\end{eqnarray*}
}
\end{remark}

\begin{table}
\begin{center}
\renewcommand{\arraystretch}{2.4}
\begin{tabular}{|c|c|}
\hline
\multicolumn{2}{|c|}{$(\, \ux \, , \, \ux\, )$}   \\ \hline
\multicolumn{2}{|c|}{$(\, r^2\, , \, r^2\, )$}    \\ \hline 
\multicolumn{2}{|c|}{$(\, \mE\, , \, \mE\, )$}  \\ \hline
$(\, \Gamma\, , \, -\, \p_{\uom}\, \uom\, )$ & $(\, -\, \p_{\uom}\, \uom\, , \Gamma\, ) $ \\ \hline
$(\, \Gamma^2\, , \, \Gamma^2 - 2(m-1)\Gamma + (m-1)^2\, )$ & $(\, \Gamma^2 - 2(m-1)\Gamma + (m-1)^2\, , \Gamma^2\, )$ \\ \hline
$(\, \dirac\, , \, \uD\, )$ & $[\, \uD\, , \, \dirac\, ]$ \\ \hline
\multicolumn{2}{|c|}{$[\, \uom\, \p_r\, , \, \uom\, \p_r\, ]$} \\ \hline
$[\, \invr\, \p_{\uom}\, , \, -\, \invr\, \p_{\uom} + (m-1)\, \invr\, \uom\, ]$ & $[\, -\, \invr\, \p_{\uom} + (m-1)\, \invr\, \uom\, , \, \invr\, \p_{\uom}\, ]$ \\ \hline
\multicolumn{2}{|c|}{$(\, \p_{\uom}^2\, , \, \p_{\uom}^2\, )$} \\ \hline
$(\, \Delta^*\, , \, {\bf Z}^*\, )$ & $(\, {\bf Z}^*\, , \, \Delta^*\, )$ \\ \hline
$(\, \Delta\, , \, {\bf Z}\, )$ & $[\, {\bf Z}\, , \, \Delta\, ]$ \\ \hline
\multicolumn{2}{|c|}{$[\, \p_r^2\, , \, \p_r^2\, ]$} \\ \hline
\multicolumn{2}{|c|}{$[\, \invux\, , \, \invux\, ]$} \\ \hline
\multicolumn{2}{|c|}{$[\, \invr\, \p_r\, , \, \invr\, \p_r\, ]$}  \\ \hline
\multicolumn{2}{|c|}{$[\, \invrsq\, , \, \invrsq\, ]$} \\ \hline
$(\, \p_{x_j}\, , \, d_j\,)$ & $[\, d_j\, , \, \p_{x_j}\, ]$ \\ \hline
\end{tabular}
\caption{Signum--pairs of operators}
\end{center}
\end{table}

\begin{table}
\begin{center}
\renewcommand{\arraystretch}{2.4}
\begin{tabular}{|c|c|}
\hline
\multicolumn{2}{|c|}{$(\, \uom \, , \, \uom\, )$} \\ \hline
$(\, r\, , \, -\, r\, )$  &  $(-\, r\, , \,  r\, )$ \\ \hline 
$[\, \p_r\, , \, -\, \p_r\, ]$ & $[-\, \p_r\, , \,  \p_r\, ]$ \\ \hline
$(\, \p_{\uom}\, , \, \uom\, \p_{\uom}\, \uom\, )$ & $(\, \uom\, \p_{\uom}\, \uom\, , \p_{\uom}\, )$ \\ \hline
$[\, \invr\, , \, -\, \invr\, ]$ & $[-\, \invr\, , \,  \invr\, ]$ \\ \hline
$[\, \invr\, \p_{\uom}\, \uom\, , \, -\, \invr\, \uom\, \p_{\uom}\, ]$ & $[-\, \invr\, \p_{\uom}\, \uom\, , \,  \invr\, \uom\, \p_{\uom}\, ]$ \\ \hline
\end{tabular}
\caption{Cross--pairs of operators}
\end{center}
\end{table}

\vspace*{8cm}

  
\section{Conclusion}
\label{conclusion}
  

In his famous and seminal book \cite{zwart} Laurent Schwartz writes on page 51:
{\em
Using co-ordinate systems other than the cartesian ones should be done with the utmost care
}
 [our translation]. And right he is! Indeed, just consider the delta-distribution $\delta(\ux)$: it is pointly supported at the origin, it is rotation invariant: $\delta(A\, \ux) = \delta(\ux), \ \forall \, A \in {\rm SO}(m)$, it is even:  $\delta(-\ux) = \delta(\ux)$ and it is homogeneous of order $(-m)$: $\delta(a \ux) = \frac{1}{|a|^m}\, \delta(\ux)$. So in a first, naive, approach, one could think of its radial derivative $\p_r\, \delta(\ux)$ as a distribution which remains pointly supported at the origin, rotation invariant, even and homogeneous of degree $(-m-1)$.  Temporarily leaving aside the even character, on the basis of the other cited characteristics the distribution $\p_r\, \delta(\ux)$ should take the following form:
$$
\p_r\, \delta(\ux) = c_0\, \p_{x_1} \delta(\ux) + \cdots + c_m\, \p_{x_m} \delta(\ux)
$$
and it becomes immediately clear that this approach to the radial derivation of the delta-distribution is impossible since all distributions appearing at the right--hand side of the above sum  are odd and not rotation invariant, whereas $\p_r\, \delta(\ux)$ is assumed to be even and rotation invariant. It could be that $\p_r\, \delta(\ux)$ is either the zero distribution or is no longer pointly supported at the origin, but both those possibilities are unacceptable. So from the start we are warned by this example that introducing spherical co-ordinates $\ux = r \uom,\, r = |\ux|,\, \uom \in \mS^{m-1}$ makes derivation of distributions in $\mR^m$ a far from trivial action, as are, in principle  ``forbidden'', actions such as multiplication by the non--analytic functions $r$ and $\omega_j, j=1,\ldots,m$. But there is more: functional analytic considerations on the space $\mcD(\mR^m)$ of compactly supported smooth test functions expressed in spherical co--ordinates, forced us to introduce a new space of continuous linear functionals on a auxiliary space of test functions showing a singularity at the origin, for which, in \cite{bsv}, we coined the term {\em signumdistributions}, bearing in mind that $\uom = \frac{\ux}{|\ux|}$ may be interpreted as the higher dimensional counterpart to the {\em signum} function on the real line. It turns out that the actions by $r$, $\uom$, $\p_r$ and $\p_{\uom}$ map a distribution to a signumdistribution and vice versa. The basic idea behind the definition of  these actions on a distribution $T \in \mcD'(\mR^m)$, is to express the resulting signumdistributions as appropriate and ``legal" actions on $T$. So, for example, we put $\langle rT,\uom \varphi \rangle = \langle r\uom T, \varphi \rangle = \langle \ux T, \varphi \rangle, \forall \varphi \in \mcD(\mR^m)$. This idea may seem to be rather simple, but it is backed up by the functional analytic considerations of Section \ref{intro}, and it paves the way for easy to handle calculus rules as established in \cite{bsv}.\\
Of the four aforementioned actions only the radial derivative $\p_r\, T$ escapes, in general, from an unambiguous definition, but leads to an equivalent class of signumdistributions instead. Still we are able to define unambiguously $\p_r\, T$ in two particular cases: (i) when the given distribution $T$ is radial, i.e. rotation invariant, and (ii) when $T = U^{\wedge}$ is the associated distribution to a given radial signumdistribution $U$, these two particular cases being quite interesting since they correspond to two families of frequently used distributions in Clifford analysis, such as the fundamental solutions of the Laplace and the Dirac operator. The spherical co-ordinates approach to the distributions in those two families will be worked out in detail in the forthcoming paper \cite{fb}.



\end{document}